\input amstex\documentstyle{amsppt}  
\pagewidth{12.5cm}\pageheight{19cm}\magnification\magstep1
\topmatter
\title Quasisplit Hecke algebras and symmetric spaces\endtitle
\author George Lusztig and David A. Vogan, Jr.\endauthor
\address{Department of Mathematics, M.I.T., Cambridge, MA 02139}\endaddress
\thanks{Supported in part by National Science Foundation grants DMS-0758262, DMS-0967272.}\endthanks
\abstract
Let $(G,K)$ be a symmetric pair over an algebraically closed field of
characteristic different from 2 and let $\sigma$ be an automorphism
with square 1 of $G$ preserving $K$. In this paper we consider the set
of pairs $({\Cal O},{\Cal L})$ where ${\Cal O}$ is a $\sigma$-stable
$K$-orbit on the flag manifold of $G$ and ${\Cal L}$ is an irreducible
$K$-equivariant local system on ${\Cal O}$ which is ``fixed" by
$\sigma$.  Given two such pairs $({\Cal O},{\Cal L})$, $({\Cal
O}',{\Cal L}')$, with ${\Cal O}'$ in the closure $\overline{\Cal O}$ 
of ${\Cal O}$, the multiplicity space of ${\Cal L}'$ in a cohomology 
sheaf of the intersection cohomology of $\overline {\Cal O}$ with 
coefficients in ${\Cal L}$ (restricted to ${\Cal O}'$) carries an 
involution induced by $\sigma$, and we are interested in computing the 
dimensions of its $+1$ and $-1$ eigenspaces. We show that this 
computation can be done in terms of a certain module structure over a 
quasisplit Hecke algebra on a space spanned by the pairs $({\Cal 
O},{\Cal L})$ as above. 
\endabstract 
\endtopmatter    
\document 
\define\tucl{\ti{\ucl}} 
\define\uph{\un\ph} 
\define\utph{\un{\tph}}

\define\tph{\ti{\ph}}

\define\ul{\un l}

\define\si{\sim}

\define\qua{\quad}

\define\lb{\linebreak}

\define\op{\oplus}

\define\part{\partial}

\define\m{\mapsto}
\define\do{\dots}

\define\sm{\smallmatrix}
\define\esm{\endsmallmatrix}
\define\sub{\subset}    

\define\T{\times}
\define\ti{\tilde}  
\define\nl{\newline}
\redefine\i{^{-1}}

\define\un{\underline} 
\define\ov{\overline}
\define\ot{\otimes} 
\define\bbq{\overline{\QQ}_l}

\define\Ad{\text{\rm Ad}}

\define\Aut{\text{\rm Aut}}

\define\tr{\text{\rm tr}}

\define\supp{\text{\rm supp}}

\define\a{\alpha}
\redefine\b{\beta}
\redefine\c{\chi}

\define\e{\epsilon}
\define\et{\eta}

\redefine\o{\omega}
\define\p{\pi}
\define\ph{\phi}

\define\r{\rho}
\define\s{\sigma}
\redefine\t{\tau}
\define\th{\theta}

\define\x{\xi}

\define\vt{\vartheta}

\define\Th{\Theta}

\define\Ph{\Phi}
\define\Ps{\Psi}

\define\kk{\bold k}

\redefine\tt{\bold t}

\define\CC{\bold C}
\define\DD{\bold D}

\define\FF{\bold F}
\define\GG{\bold G}
\define\HH{\bold H}

\define\KK{\bold K}

\define\NN{\bold N}

\define\QQ{\bold Q}
\define\RR{\bold R}
\define\SS{\bold S}

\define\ZZ{\bold Z}

\define\ca{\Cal A}
\define\cb{\Cal B}
\define\cc{\Cal C}

\define\cf{\Cal F}
\define\cg{\Cal G}

\define\cl{\Cal L}

\define\co{\Cal O}

\define\cs{\Cal S}
\define\ct{\Cal T}

\define\cv{\Cal V}

\define\fg{\frak g}

\define\fA{\frak A}

\define\fC{\frak C}
\define\fD{\frak D}

\define\fK{\frak K}

\define\fO{\frak O}

\define\fS{\frak S}

\define\sh{\sharp}

\define\bS{\bar S}

\define\ucl{\un\cl}

\define\Un{ALTV}
\define\IW{Iw}
\define\KL{KL}
\define\MA{Ma}
\define\LC{L1}  
\define\CS{L2}  
\define\QG{L3}
\define\DG{L4}
\define\LV{LV1}
\define\LVV{LV2}
\define\BBD{BBD}
\hyphenation{con-struc-ti-ble}
\head Introduction\endhead

Suppose $\GG$ is a complex connected reductive algebraic 
group. \'Elie Cartan showed that the equivalence classes of real forms
of $\GG$ are in one-to-one correspondence with equivalence classes of 
automorphisms $\th\colon \GG @>>>\GG$ such that $\th^2=1$. Suppose $\th$ is
such an automorphism, so that $\KK = \GG^\th$ is a (possibly disconnected)
complex reductive subgroup of $\GG$. Cartan's correspondence has 
the property that $\GG(\RR)$ contains the (unique) compact real form
$\KK(\RR)$ as maximal compact subgroup. For example, if $\GG =
GL(n,\CC)$, then the unique compact real form of $\GG$ is $U(n)$; this
corresponds to the trivial automorphism $\th$. The automorphism
$\th(g) = {}^tg^{-1}$ has fixed points the complex orthogonal group
$O(n,\CC)$; the corresponding real form is $GL(n,\RR)$, which contains
$O(n,\RR)$ as maximal compact subgroup.

Problems of harmonic analysis on $\GG(\RR)$---which arise for example 
in the theory of automorphic forms, or in the study or 
$\GG(\RR)$-invariant differential equations---lead to the study of 
infinite-dimensional representations of $\GG(\RR)$. Harish-Chandra 
showed in the 1950s that these difficult analytic objects could be 
related to algebraic ones: ``$(\fg,\KK)$-modules,'' which are 
simultaneously algebraic representations of the algebraic group $\KK$
and representations of the Lie algebra $\fg$. Finally the localization
theorem of Beilinson and Bernstein related $(\fg,\KK)$-modules to
algebraic geometry: $\KK$-equivariant perverse sheaves on the flag
variety of $\GG$.  To be slightly more explicit, certain Euler
characteristics of the local cohomology of these perverse sheaves are
coefficients in Weyl-type formulas for Harish-Chandra's distribution
characters of the corresponding $\GG(\RR)$ representations.

The paper \cite{\LV} provides a method to compute the dimensions of
these local cohomology groups for $\KK$-equivariant intersection
cohomology complexes; and, consequently, character formulas for
irreducible representations of $\GG(\RR)$. We recall the idea. Very
general arguments from arithmetic geometry (see for example
\cite{\BBD, \S6.1}) allow us to replace the complex constructible
sheaves related to $\GG$ and $\KK$ by corresponding $l$-adic
constructible sheaves on related groups over a field
$\overline \FF_q$. (There is a lot of flexibility in the prime power
$q$; it can be taken to be any sufficiently large power of almost any
prime.) Precisely, {\it the dimensions of stalks of complex
intersection homology sheaves that we want are equal to the dimensions
of stalks of $l$-adic intersection homology sheaves in characteristic
$p$}.  It is this last setting ($l$-adic sheaves in odd characteristic
$p\ne l$) that is studied in \cite{\LV}.  What is gained by the
translation is that one can work with not just dimensions, but with
characteristic polynomials of Frobenius automorphisms; and for these,
there are deep and powerful tools available.

In the present paper, we are interested in a second automorphism $\s$
of the complex reductive group $\GG$, also assumed to satisfy $\s^2 =
1$ and to commute with $\th$. (Allowing $\s$ of other finite orders
requires no essentially new ideas, but complicates significantly
the detailed computations of \S 7. Because it suffices for the
applications described next, we consider only the case $\s^2 = 1$.)

Just as in \cite{\DG}, it is often
convenient to think of $\s$ as defining a {\it disconnected} reductive
group $\hat\GG = \GG \rtimes \{1,\sigma\}$; the assumption that $\s$
commutes with $\th$ means that this disconnected group is also defined
over $\RR$. It is shown in \cite{\Un} that a precise understanding of
the (infinite-dimensional) representation theory of $\hat\GG(\RR)$
(for an appropriately chosen $\sigma$) leads to an algorithm for
determining the {\it unitary} irreducible representations of
$\GG(\RR)$.

Clifford theory says that the representation theory of $\hat\GG(\RR)$
is very close to that of $\GG(\RR)$. The main point is to understand
the action of $\s$ on irreducible representations of $\GG(\RR)$. The
translations of Harish-Chandra (to algebra) and Beilinson-Bernstein
(to geometry) remain valid: the conclusion is that irreducible
representations of $\GG(\RR)$ may be described in terms of
$\hat\KK$-equivariant perverse sheaves on the flag manifold for $\GG$;
that is, in terms of the action of $\s$ on stalks of local cohomology
of perverse sheaves. Just as in the classical case, there is a
straightforward translation of the problem to any finite
characteristic $p\ne 2$.  {\it The problem we address is therefore
(essentially) computation of the trace of $\s$ acting on stalks of the
intersection cohomology sheaves considered in \cite{\LV}.} 
Appropriate sums of these traces amount to coefficients in character
formulas for irreducible representations of the disconnected reductive
group $\hat\GG(\RR)$, and consequently play a role (explained in
\cite{\Un}) in determining the unitary representations of $\GG(\RR)$.

The main idea in \cite{\LV} was to make an appropriate Grothendieck
group of a category of equivariant sheaves into a module for the
Iwahori Hecke algebra of the Weyl group $W$ of $\GG$; to calculate
this module action explicitly in a basis of sheaves supported on
single orbits; and to establish a relationship between the action and
Verdier duality on complexes of sheaves.  Combining all this
information gave an algorithm for computing the intersection homology
sheaves.

In the setting of the present paper, the action we need is of a
smaller, unequal parameter Hecke algebra: one related to the fixed
points of a natural action of $\s$ on $W$.  The group $W^\s$ appears
as the Weyl group of the maximally split torus in a quasisplit finite
Chevalley group (having the same root datum as $\GG$). The action is
introduced in 4.2(g). With this action 
in hand, the strategy of proof is roughly the same as in \cite{\LV};
the key statements relating Verdier duality and the Hecke algebra
action are in 4.8. 

\subhead 0.1\endsubhead
We turn now to a more precise description of our results. Let $G$ be a
connected reductive algebraic group over an algebraic 
closure $\kk$ of the finite field $\FF_p$ with 
$p$ elements ($p$ is a prime number $\ne2$). We assume that we are
given a closed subgroup $K$ of $G$ such that 
$(G,K)$ is a symmetric pair (thus, we are given an automorphism
$\th\colon G@>>>G$ such that $\th^2=1$ and $K$ has  
finite index in $G^\th$, the fixed point set of $\th$) and an
automorphism $\s\colon G@>>>G$ such that $\s^2=1$ and  
$\s\th=\th\s$, $\s(K)=K$. Sometimes it is convenient to think of 
$\sigma$ as defining a disconnected reductive group 
$$\hat G = G \rtimes \{1,\sigma\} \supset \hat K = K \times 
\{1,\sigma\}$$ as in \cite{\DG}. 
Let $\cb$ be the variety of Borel subgroups of $G$. Then $K$ acts on
$\cb$ by  
conjugation with finitely many orbits that form a set $E$. 
 
We fix a prime number $l$ such that $l\ne p$. Let $\cc_0$ be the
category whose objects are the constructible  
$K$-equivariant $\bbq$-sheaves on $\cb$; the morphisms in $\cc_0$ are
assumed to be compatible with the 
$K$-equivariant structures. Let $\fD$ be a set of representatives for
the isomorphism classes of objects  
$\cl\in\cc_0$ such that for some $K$-orbit $\co$ in $\cb$, $\cl|_\co$
is a $\bbq$-local system of rank $1$ and  
$\cl|_{\cb-\co}=0$. Note that $\co$ is uniquely determined by $\cl$
and is denoted by $[\cl]$; we shall write  
$\ucl$ instead of $\cl|_{[\cl]}$. Note that $\fD$ is a finite set. For
any $\cs\in\cc_0$ and any $\co\in E$ we have a canonical decomposition 
$$\cs|_\co=\bigoplus_{\cl\in\fD;[\cl]=\co}V_\cl(\cs)\ot\ucl$$
(as $K$-equivariant local systems over $\co$) where $V_\cl(\cs)$ are
finite-dimensional $\bbq$-vector spaces. 

For any complex $L$ of constructible $\bbq$-sheaves on an algebraic
variety we denote by $L^i$ the  $i$-th cohomology sheaf of $L$. 

For $\cl\in\fD$ let $\cl^\sh$ be the intersection cohomology complex
of the closure $\ov{[\cl]}$ of $[\cl]$  
with coefficients in $\ucl$, extended by $0$ on $\cb-\ov{[\cl]}$. For
any $i\in\ZZ$, $\cl^{\sh i}$ is naturally  
an object of $\cc_0$; it is zero for all but finitely many $i$ and is zero unless $i\in2\NN$, see \cite{\LV}. For
any $\cl,\cl'\in\fD$ we set 
$$P_{\cl',\cl}=\sum_{h\in\NN}\dim V_{\cl'}(\cl^{\sh 2h})u^h\in\NN[u]$$
($u$ is an indeterminate). Now the polynomials $P_{\cl',\cl}$ were studied in \cite{\LV}; they are of interest
for the representation theory of real reductive groups.

In this paper we will consider a variant of these polynomials in which $\s$ plays a role. If $B\in\cb$ then
$\s(B)\in\cb$ and $B\m\s(B)$ defines an involution of $\cb$ denoted by $\s$. This induces a permutation of $E$
whose fixed point set is denoted by $E^\s$. Let $\fD^\s$ be the set of all $\cl\in\fD$ such that $\s^*\cl$ is 
isomorphic to $\cl$ in $\cc_0$; if $\cl\in\fD^\s$, then $[\cl]\in E^\s$. For any $\cl\in\fD^\s$ let $N_\cl$ be 
the set of isomorphisms $\a\colon \s^*\cl@>\si>>\cl$ in $\cc_0$ such
that $\a\s^*(\a)\colon \cl@>>>\cl$ is the  
identity. Note that $|N_\cl|=2$. (In the language of disconnected 
groups, the choice of $\a \in N_\cl$ is the 
same as the choice of a $\hat K$-equivariant structure on $\cl$ 
extending the given $K$-equivariant structure.) For each
$\cl\in\fD^\s$ we select $\a^\cl\in N_\cl$; let  
$\a^{\cl\sh}\colon \s^*\cl^\sh@>\si>>\cl^\sh$ be the canonical
extension of $\a^\cl$; this induces for any $h\in\NN$ an 
isomorphism $\a^{\cl\sh2h}\colon
\s^*\cl^{\sh2h}@>\si>>\cl^{\sh2h}$. If $\co\in E^\s$ then  
$\a^{\cl\sh2h}|_\co$ can be viewed as an isomorphism  
$$\bigoplus_{\cl'\in\fD;[\cl']=\co}V_{\cl'}(\cl^{\sh2h})\ot\s^*\ucl'@>>>
 \bigoplus_{\cl'\in\fD;[\cl']=\co}V_{\cl'}(\cl^{\sh2h})\ot\ucl'$$
which maps the $\cl'$-summand on the left to the $\cl'$-summand on the right (if $\cl'\in\fD^\s$) according to an
isomorphism of the form $\a^{\cl\sh2h;\cl'}\ot\a^{\cl'}$ where 
$$\a^{\cl\sh2h;\cl'}\colon
V_{\cl'}(\cl^{\sh2h})@>>>V_{\cl'}(\cl^{\sh2h})$$
is a vector space isomorphism whose square is
$1$. For any $h\in\NN$ we set
$$P^\s_{\cl',\cl;h}=\tr(\a^{\cl\sh2h;\cl'}\colon V_{\cl'}(\cl^{\sh2h})@>>>V_{\cl'}(\cl^{\sh2h}))\in\ZZ.$$
We also set 
$$P^\s_{\cl',\cl}=\sum_{h\in\NN}P^\s_{\cl',\cl;h}u^h\in\ZZ[u].$$
Thus $P^\s_{\cl',\cl}$ is defined in the same way as $P_{\cl',\cl}$, but using trace instead of dimension. Note 
that $P^\s_{\cl',\cl}$ depend on the choice of $\a^\cl$ for each $\cl\in\fD^\s$; another choice can change
$P^\s_{\cl',\cl}$ to its negative. The polynomials $P^\s_{\cl',\cl}$ are expected to be of interest for the 
theory of {\it unitary} representations of real reductive groups. 

A special case of these polynomials was considered in \cite{\LC} where
$G$ was replaced by $G\T G$, $\th$ was  
replaced by the map $(g,g')\m(g',g)$ and $\s$ was replaced by the map
$(g,g')\m(\s(g),\s(g'))$ (see 0.2). The  
polynomials introduced in \cite{\LV} can be also viewed as a special
case of the polynomials in the present  
paper; they correspond to the case where $\s=1$. (The polynomials in
\cite{\LC} and those in \cite{\LV}  
generalize those in \cite{\KL} in different directions.) More  
recently, another special case of these polynomials was considered in
\cite{\LVV}, where $G$ was replaced by 
$G\T G$, $\th$ was replaced by the map $(g,g')\m(g',g)$ and $\s$ was
replaced by the map $(g,g')\m(\s(g'),\s(g))$. 
(In each of the three special cases above there was a canonical choice
for $\a^\cl$.) 

One of the main results of this paper is the construction of an action
of a (quasisplit) Hecke algebra on a 
module spanned by the elements of $\fD^\s$. To do this we use, in
addition to the techniques of \cite{\LV}, an  
idea from the geometric construction \cite{\QG, Ch. 12} of the plus
part of a universal quantized enveloping  
algebra of nonsimplylaced type as a quotient of a Grothendieck group
associated to a periodic functor on a 
category. (This idea was also used in \cite{\LVV}.)  

We also show how the polynomials $P^\s_{\cl',\cl}$ can be characterized in terms of a certain bar operator on 
this Hecke-module (see Theorem 5.2). This reduces the problem of explicitly computing the $P^\s_{\cl',\cl}$ to 
the problem of explicitly computing the bar operator. 

In \S 7 we compute explicitly the action of the generators of the
Hecke algebra ${\bold H}$ on the basis of $M$ parametrized by ${\Cal
D}^\sigma$. The computation is facilitated by a model
(described in \S 6) for (the specialization to an odd power of $q$) of
$M$ in terms of certain functions on the rational points of ${\Cal
B}$.  Using this model, the computation of the action of generators
can in the most complicated cases be reduced to groups locally
isomorphic to $SL(2)$, $SL(2) \times SL(2)$, and $SL(3)$.  A
description of these special cases is sketched in \S 9.

In \S 8, we explain how the formulas of \S 7 lead to a recursive 
algorithm for calculating the bar operator on $M$, or (equivalently) 
the polynomials $P^s_{\cl',\cl}$. In the setting of \cite{\KL}, what 
makes the recursion work is that the identity element---on which the 
bar operation acts trivially---is a generator for the Hecke algebra as
a module over itself.  Similarly, in \cite{\LVV}, the Hecke module 
corresponding to twisted involutions is generated over the quotient 
field by elements corresponding to orbits of minimal dimension. In the
present case, as in \cite{\LV}, the Hecke module $M$ need not be 
generated by local systems on orbits of minimal dimension, or even by 
local systems fixed by bar. The resolution we give here is a little 
different.  It is based on unpublished work of the second author with 
Jeffrey Adams and Peter Trapa in the setting of \cite{\LV}. 

\subhead 0.2\endsubhead 
We now consider the following special case. We replace $G$ by $G\T G$
(hence $\cb$ is replaced by $\cb\T\cb$). We 
replace $\th$ by the map $(g,g')\m(g',g)$ (hence $K$ is replaced by
the diagonal $G$ in $G\T G$). 
We replace $\s$ by $\s\T\s$. Then $\cc_0$ becomes $\fC_0$ (the objects
of $\fC_0$ are the  
constructible $G$-equivariant $\bbq$-sheaves on $\cb\T\cb$). Let $W$
be the set of orbits of the $G$-action on 
$\cb\T\cb$. This is naturally a finite Coxeter group (the Weyl group)
with a standard length function  
$\ul\colon W@>>>\NN$. Let $S$ be the standard set of generators of $W$. For
$w\in W$ we write $\fO_w$ for the  
corresponding $G$-orbit in $\cb\T\cb$. We can identify $\fD$ with $W$
(to $w\in W$ corresponds the object $\SS_w\in\fC_0$ such that
$\SS_w|_{\fO_w}=\bbq$, 
$\SS_w|_{(\cb\T\cb)-\fO_w}=0$).  
In our case the subset $\fD^\s$ of $\fD$ becomes the subgroup
$W^\s=\{w\in W;\s(w)=w\}$ of $W$ (note that $\s$ 
acts naturally on $W$). Note that for $w\in W^\s$ there is a canonical
choice for the isomorphism 
$\a^{\SS_w}\colon (\s\T\s)^*\SS_w@>>>\SS_w$ namely the one inducing the
identity map on the stalks over $\fO_w$ (all 
these stalks are $\bbq$).  
For $w\in W^\s$ let $\SS_w^\sh$ be the intersection cohomology
complex of the closure $\ov{\fO_w}$ of $\fO_w$ with coefficients in
$\bbq$, extended by $0$ on  
$(\cb\T\cb)-\ov{\fO_w}$. We have $\SS_w^{\sh i}=0$ unless $i\in2\NN$
(\cite{\KL}).  

Now let $y,w\in W^\s$ be such that $y\le w$ where $\le$ is the Bruhat
order on $W$. Then $\fO_y$ is contained in $\ov{\fO_w}$, the closure
of $\fO_w$. The $(2h)$-th cohomology sheaf of $\SS_w^\sh$ restricted
to $\fO_y$ is a local system whose stalks can be identified with a
single $\bbq$-vector 
space $V_{y,w;2h}$ which carries a canonical involution induced by
$\s$. Let $n_{y,w;2h}$ be the trace of this  involution. Let
$P^\s_{y,w}=\sum_{h\in\NN}n_{y,w;2h}u^h\in\ZZ[u]$. This polynomial
(which is a special case of the polynomials $P^\s_{\cl',\cl}$ in 0.1)
was considered in \cite{\LC, (8.1)} where it was stated 
without proof that it can be explicitly computed by an algorithm
involving the quasisplit Iwahori-Hecke algebra  
$H$ associated to the Coxeter group $W^\s$ and the restriction of
$\ul$ to $W^\s$. (In the case where $\s=1$ the  
trace becomes dimension and the statement of \cite{\LC, (8.1)} reduced
to a result in \cite{\KL}.) A proof of the 
statement for general $\s$ was given in \cite{\DG} where it was also
shown how the polynomials $P^\s_{y,w}$ enter 
in the study of unipotent representations of disconnected reductive
groups over a finite field. The methods of  
this paper give another proof of this statement and at the same time
give a geometric construction of the algebra 
$H$ which is one of the ingredients in our construction of the
$H$-module in 0.1. 

\subhead 0.3\endsubhead
{\it Notation.} For an algebraic group $H$ we denote by $H^0$ the identity component of $H$. If $X$ is a set and 
$f\colon X@>>>X$ is a map, we define $$X^f=\{x\in X;f(x)=x\}.$$ 
For any $q$ (a power of $p$) let $\FF_q$ be the subfield 
of $\kk$ such that $|\FF_q|=q$. If $B\in\cb$ we set $K_B=K\cap B$;
then $K_B/K_B^0$ is an elementary  
abelian $2$-group. If $\cs\in\cc_0$ and $B\in\cb$, let $\cs_B$ be the
stalk of $\cs$ at $B$; if $x\in K$, we  
denote by $\ct^\cs_x\colon \cs_B@>>>\cs_{xBx\i}$ the linear isomorphism given by the $K$-equivariant structure of 
$\cs$. For any complex $L$ of constructible $\bbq$-sheaves on an
algebraic variety we denote by $DL$ the Verdier  
dual of $L$. If $\cs\in\cc_0$ then $D^i\cs\colon =(D\cs)^i$ is naturally an object of $\cc_0$; it is zero for all but 
finitely many $i$. 

\head 1. The category $\cc'$\endhead
\subhead 1.1\endsubhead
In this section we review (and slightly strengthen) some results of \cite{\LV}.

We can find (and will fix) a morphism $\ph\colon G@>>>G$ which is the Frobenius map for a split $\FF_q$-rational 
structure on $G$ ($q$ is a sufficiently divisible power of $p$) such that (denoting the map $\cb@>>>\cb$, 
$B\m\ph(B)$, again by $\ph$) the following hold.

(i) $\ph\th=\th\ph$ and $\ph(K)=K$;

(ii) any $K$-orbit on $\cb$ meets $\cb^\ph$ (hence is $\ph$-stable);

(iii) for any $\cl\in\fD$ we have $\ph^*\cl\cong\cl$ in $\cc_0$;

(iv) we have $\ph\s=\s\ph$.
\nl
Note that condition (iv) will not play any role in this section. On the other hand in 2.1 we will add another
requirement to (i)-(iv) above.

Now, if $\cs\in\cc_0$ then $\ph^*\cs$ is naturally an object of $\cc_0$ such that for any $k\in K,B\in\cb$, 
$\ct^{\ph^*\cs}_k\colon (\ph^*\cs)_B@>>>(\ph^*\cs)_{kBk\i}$ is the same as 
$\ct^\cs_{\ph(k)}\colon \cs_{\ph(B)}@>>>\cs_{\ph(kBk\i)}$. (We use that $\ph(K)=K$.)

\proclaim{Lemma 1.2}Let $\cl\in\fD$ and let $\co=[\cl]$.

(a) The stalks $\ucl^{\ot2}_B$ for various $B\in\co$ can be canonically identified with a single $\bbq$-vector 
space $\cv$ independent of $B$.

(b) There exists an isomorphism $t\colon \ph^*\cl@>\si>>\cl$ in $\cc_0$ such that for any $B\in\cb^\ph$, the induced 
map $t_B\colon \cl_B@>>>\cl_B$ is the identity map or $(-1)$ times the identity map. Moreover if 
$t'\colon \ph^*\cl@>\si>>\cl$ is an isomorphism with the same properties then $t'=t$ or $t'=-t$.
\endproclaim
We prove (a). Let $B,B'\in\co$. If $k,k'\in K$ are such that $kBk\i=B'$, $k'Bk'{}\i=B'$ then
$\ct_{k'}^\cl=\ct_k^\cl\ct_{k_0}^\cl$ where $k_0=k'k\i\in K_B$ and $\ct_{k_0}^\cl\colon \ucl_B@>>>\ucl_B$ is in the 
image of a homomorphism $K_B@>>>\Aut(\ucl_B)$ whose kernel contains $K^0_B$. Since $K_B/K^0_B$ is an 
elementary abelian $2$-group and $\dim\ucl_B=1$ we see that $\ct_{k_0}^\cl=\pm1$. Thus, 
$\ct_{k'}^\cl=\pm\ct_k^\cl$. It follows that $(\ct_k^\cl)^{\ot2}\colon \ucl^{\ot2}_B@>>>\ucl^{\ot2}_{B'}$ is 
independent of the choice of $k\in K$ such that $kBk\i=B'$. Thus we have a canonical isomorphism 
$\ct_{B,B'}\colon \ucl^{\ot2}_B@>\si>>\ucl^{\ot2}_{B'}$. This has an obvious transitivity property; (a) follows.

We prove (b). By 1.1(iii) we can find an isomorphism of $K$-equivariant local systems $t\colon \ph^*\ucl@>>>\ucl$. We
can assume that for some $B_0\in\co^\ph$, $t_{B_0}\colon \ucl_{B_0}@>>>\ucl_{B_0}$ is the identity map (see 1.1(ii)). 
Now $t$ induces an isomorphism $t^{\ot2}\colon \ph^*\ucl^{\ot2}@>>>\ucl^{\ot2}$. For any $B\in\co$, the induced map on 
stalks $t^{\ot2}_B\colon \ucl^{\ot2}_{\ph(B)}@>>>\ucl^{\ot2}_B$ can be viewed as a linear isomorphism $\cv@>>>\cv$ (see
(a)) which is independent of $B$ and is multiplication by a scalar $c\in\bbq^*$. Taking $B=B_0$ we see that 
$c=1$. Thus for any $B\in\co^\ph$, $t^{\ot2}_B\colon \ucl^{\ot2}_B@>>>\ucl^{\ot2}_B$ is the identity map, hence
$t_B\colon \cl_B@>>>\cl_B$ is multiplication by a scalar whose square is $1$. This proves the existence part of (b). 
Now let $t'$ be as in (b). Since $\ucl$ is irreducible as a $K$-equivariant local system, there exists 
$b\in\bbq^*$ such that $t'=bt$. By assumption, if $B\in\co^\ph$, $t_B\colon \cl_B@>>>\cl_B$ and $t'_B\colon \cl_B@>>>\cl_B$ 
are multiplication by $\pm1$; hence $t'_B=\e t_B$ where $\e=\pm1$. It follows that $b=\e$. This proves (b). The 
lemma is proved.

\subhead 1.3\endsubhead
Let $\cc_1$ be the category whose objects are pairs $(\cs,t)$ where $\cs\in\cc_0$ and $t$ is an isomorphism 
$\ph^*\cs@>\si>>\cs$ in $\cc_0$ such that for any $B\in\cb^\ph$, the
eigenvalues of the induced linear map  
$t_B\colon \cs_B=\cs_{\ph(B)}@>>>\cs_B$ are of the form $\pm q^e (e\in\ZZ)$. A morphism between two objects 
$(\cs,t),(\cs',t')$ of $\cc_1$ is a morphism $\cs@>e>>\cs'$ in $\cc_0$
such that the diagram 
$$\CD \ph^*\cs@>t>>\cs\\
@V\ph^*e VV     @V e VV\\
\ph^*\cs'@>t'>>\cs'\endCD$$
is commutative. For example, if $\cl\in\fD,k\in\ZZ$ and $t$ is as in
1.2(b), then $(\cl,q^kt)$ is an object of $\cc_1$. 

\proclaim{Lemma 1.4}Let $\cl\in\fD$ and let $t$ be as in 1.2(b). Let
$t^\sh\colon \ph^*\cl^\sh@>>>\cl^\sh$ be the  
isomorphism which extends $t\colon \ph^*\cl@>>>\cl$. For any $h\in\NN$ let 
$t^{\sh2h}\colon \ph^*(\cl^{\sh 2h})@>>>\cl^{\sh 2h}$ be the
isomorphism induced by $t^\sh$. Then for any  
$B\in\cb^\ph$, any eigenvalue of $t^{\sh2h}$ on $\cl^{\sh 2h}_B$ is
equal to $q^h$ or to $-q^h$. In particular we 
have $(\cl^{\sh 2h},t^{\sh2h})\in\cc_1$. 
\endproclaim 
The weaker statement that any eigenvalue of $t^{\sh2h}$ on
$\cl^{\sh2h}_B$ is equal to $q^h$ times a root of $1$ 
is given in \cite{\LV,4.10}. The main reason that {\it loc.cit.}
gives only the weaker statement is that the 
statement 1.2(b) was not available there. But once 1.2(b) is known we
can essentially repeat the arguments in  
{\it loc.cit.} and obtain the desired result. 

\proclaim{Lemma 1.5} Let $(\cs,t)\in\cc_1$ and let $i\in\ZZ$. Let
$t^{(i)}\colon \ph^*(D^i\cs)@>>>D^i\cs$ be the  
inverse of the isomorphism $D^i\cs@>>>\ph^*(D^i\cs)$ induced by
$D$. We have $(D^i\cs,t^{(i)})\in\cc_1$. 
\endproclaim
It suffices to show this assuming that $\cs=\cl\in\fD$ and $t$ is as
in 1.2(b). Let $\co=[\cl]$. We can assume  
that the result is known when $\cl$ is replaced by $\cl'$ where
$[\cl']\sub\ov{\co}-\co$. Let  
$t^\sh\colon \ph^*\cl^\sh@>>>\cl^\sh$ be as in 1.4. Let $L$ be
$\cl^\sh|_{\ov{\co}-\co}$ extended by $0$ on  
$\cb-(\ov{\co}-\co)$. Now $t^\sh$ induces an isomorphism $d\colon
\ph^*L@>\si>>L$. Let $d'\colon \ph^*DL@>\si>>DL$ be the  
inverse of the isomorphism $DL@>>>\ph^*(DL)$ induced by $d$; this
induces isomorphisms   
$d'_i\colon \ph^*(D^iL)@>\si>>D^iL$ for $i\in\ZZ$. Since
$\supp(L)\sub\ov{\co}-\co$ we see using the induction hypothesis that  

(a) $((DL)^h,d'_h)\in\cc_1$ for $h\in\ZZ$. 
\nl 
Now $t^\sh$ induces an isomorphism $D\cl^\sh@>>>\ph^*D\cl^\sh$ whose
inverse is an isomorphism 
$j\colon \ph^*D\cl^\sh@>\si>>D\cl^\sh$. This induces for any $i\in\ZZ$
an isomorphism 
$j^i\colon \ph^*(D^i\cl^\sh)@>\si>>D^i\cl^\sh$. We can identify
$D\cl^\sh= \cl^\sh[2m]$ for some $m\in\ZZ$ in such 
a way that $j^i$ becomes $q^{2m'}t^{\sh i+2m}$ for some
$m'\in\ZZ$. Using Lemma 1.4 we deduce that 

(b) $(D^i\cl^\sh,j^i)\in\cc_1$. 
\nl
Using (b), (a) and the long exact sequence of cohomology sheaves
associated to the exact triangle consisting of $D\cl,D(\cl^\sh),DL$
(which is obtained from the exact triangle consisting of
$L,\cl^\sh,\cl$) we deduce that  $(D^i\cl,t^{(i)})\in\cc_1$ for any
$i\in\ZZ$. This completes the inductive proof. 

\subhead 1.6\endsubhead 
Let $\cl\in\fD$. We define $\tt^\cl\colon \ph^{2*}\cl@>\si>>\cl$ as the composition 
$\ph^{2*}\cl@>\ph^*t>>\ph^*\cl@>t>>\cl$ where $t\colon \ph^*\cl@>\si>>\cl$ is as in 1.2(b). Note that by 1.2(b), 
$\tt^\cl$ is independent of the choice of $t$. 

Let $\cc'$ be the category whose objects are pairs $(\cs,\Ps)$ where
$\cs\in\cc_0$ and $\Ps$ is an isomorphism  
$\ph^{2*}\cs@>\si>>\cs$ in $\cc_0$ such that the equivalent conditions
(i), (ii) below are satisfied (for any  $\co\in E$, we identify
$\Ps|_\co$ with an isomorphism 
$$ \bigoplus_{\cl\in\fD;[\cl]=\co} V_\cl(\cs)\ot\ph^{2*}
\cl@>>> \bigoplus_{\cl\in\fD;[\cl]=\co} V_\cl(\cs)\ot\cl$$
of the form $\op_{\cl\in\fD;[\cl]=\co}\Ps^\cl\ot\tt^\cl$ where
$\Ps^\cl\colon V_\cl(\cs)@>\si>>V_\cl(\cs)$ is a vector space isomorphism).

(i) For any $B\in\cb^\ph$, the eigenvalues of the induced map $\Ps_B\colon \cs_B=\cs_{\ph^2B}@>>>\cs_B$ are of the form
$q^{2e} (e\in\ZZ)$. (Note that $\cb^\ph\sub\cb^{\ph^2}$.)

(ii) For any $\cl\in\fD$ the eigenvalues of $\Ps^\cl\colon V_\cl(\cs)@>\si>>V_\cl(\cs)$ are of the form 
$q^{2e} (e\in\ZZ)$.
\nl
A morphism between two objects $(\cs,\Ps)$, $(\cs',\Ps')$ of $\cc'$ is
a morphism $\cs@>e>>\cs'$ in $\cc_0$ such that the diagram
$$\CD\ph^{2*}\cs@>\Ps>>\cs\\ 
@V\ph^{2*}e VV     @V e VV\\
\ph^{2*}\cs'@>\Ps'>>\cs'\endCD$$ 
is commutative. 

(a)  Let $\cl\in\fD$. Then $(\cl,\Ps)\in\cc'$ for a unique $\Ps$. We
have $\Ps=\tt^\cl$. 
\nl 
This is immediate since by 1.2(b), for any $B\in\cb^\ph$,
$\tt^\cl_B\colon \cl_B@>>>\cl_B$ is the identity map.   

(b) Let $\tt^{\cl\sh}\colon \ph^{2*}\cl^\sh@>>>\cl^\sh$ be the
isomorphism which extends $\tt^\cl\colon \ph^{2*}\cl@>>>\cl$.  
For any $h\in\NN$ let $\tt^{\cl\sh2h}\colon
\ph^{2*}(\cl^{\sh2h})@>>>\cl^{\sh2h}$ be the isomorphism induced by  
$\tt^{\cl\sh}$. Then for any $B\in\cb^\ph$, any eigenvalue of
$\tt^{\cl\sh2h}$ on $\cl^{\sh2h}_B$ is equal to  
$q^{2h}$. (This follows from Lemma 1.4 since $\tt^{\cl\sh2h}\colon
\cl^{\sh2h}_B@>>>\cl^{\sh2h}_B$ is the square of  
$t^{\sh2h}\colon \cl^{\sh2h}_B@>>>\cl^{\sh2h}_B$.) Hence for any
$\cl'\in\fD$, any eigenvalue of $\tt^{\cl\sh2h;\cl'}\colon
V_{\cl'}(\cl^{\sh2h})@>>>V_{\cl'}(\cl^{\sh2h})$ is equal to
$q^{2h}$. In particular we have
$(\cl^{\sh2h},\tt^{\cl\sh2h})\in\cc'$. 

We have the following variant of Lemma 1.5. 

(c) Let $(\cs,\Ps)\in\cc'$ and let $i\in\ZZ$. Let $\Ps^{(i)}\colon \ph^{2*}(D^i\cs)@>>>D^i\cs$ be the inverse of the 
isomorphism $D^i\cs@>>>\ph^{2*}(D^i\cs)$ induced by $D$. We have $(D^i\cs,\Ps^{(i)})\in\cc'$.
\nl
It suffices to show this assuming that $\cs=\cl\in\fD$ and $\Ps=\tt^\cl$. In this case, for any $B\in\cb^\ph$, 
the linear map $(D^i\cl)_B@>>>(D^i\cl)_B$ induced by $\Ps^{(i)}$ is the square of the linear map 
$(D^i\cl)_B@>>>(D^i\cl)_B$ induced by $t^{(i)}$ where $t\colon \ph^*\cl@>>>\cl$ is as in 1.2(b). Hence the result 
follows from Lemma 1.5.

\subhead 1.7\endsubhead
Let $y\in W$. We consider the diagram $\cb@<\p_1<<\fO_y@>\p_2>>\cb$ where $\p_1(B,B')=B,\p_2(B,B')=B'$. For 
$(\cs,\Ps)\in\cc'$ and $i\in\ZZ$ we have naturally $(\p_{1!}\p_2^*\cs)^i\in\cc_0$ (we use that $K$ acts 
naturally on $\fO_y$ so that $\p_1,\p_2$ are $K$-equivariant); we denote by 
$\Ps_i\colon \ph^{2*}(\p_{1!}\p_2^*\cs)^i@>>>(\p_{1!}\p_2^*\cs)^i$ the isomorphism induced by $\Ps$ (we use that $\fO_y$ 
has a natural $\FF_q$-structure such that $\p_1,\p_2$ are defined over $\FF_q$).

\proclaim{Lemma 1.8} We have $((\p_{1!}\p_2^*\cs)^i,\Ps_i)\in\cc'$.
\endproclaim
We argue by induction on $\ul(y)$. If $y=1$ there is nothing to prove. Assume now that $y\in S$. By a standard 
argument we can assume that $(\cs,\Ps)=(\cl,\tt^\cl)$ (with $\cl\in\fD$). We denote by 
$t_i\colon \ph^*(\p_{1!}\p_2^*\cs)^i@>>>(\p_{1!}\p_2^*\cs)^i$ 
the isomorphism induced by $t\colon \ph^*\cl@>>>\cl$ as in 1.2(b). It is enough to show that 
$((\p_{1!}\p_2^*\cl)^i,t_i)\in\cc_1$. This is implicit in the proof of \cite{\LV, Lemma 3.5}. Next we assume that 
$\ul(y)\ge2$ and that the result is known for elements of length $<\ul(y)$. We can find $s\in S$ and $y'\in W$ 
such that $y=sy'$, $\ul(y)=\ul(y')+1$. We have a diagram
$$\matrix
&&\hskip -4pt\fO_s&&&& \hskip -4pt\fO_{y'}&& \\ \vspace {-.5ex}
& \hskip -8pt\raise.6ex\hbox{${}^{\p'_1}$}\hskip -8pt\swarrow &&\hskip
-8pt \searrow\hskip
-5pt\raise.6ex\hbox{${}^{\p''_1}$} &&
\hskip -6pt\raise.6ex\hbox{${}^{\p'_2}$}\hskip -6pt\swarrow &&  
\hskip -8pt\searrow\hskip -5pt\raise.6ex\hbox{${}^{\p''_2}$} &\\ \vspace{.5ex} 
\cb&&&& \hskip -8pt \cb&&&&\hskip -16pt\cb 
\endmatrix
$$
where $\p'_1(B_1,B_2)=B_1$, $\p''_1(B_1,B_2)=B_2$, $\p'_2(B_1,B_2)=B_1$, $\p'_2(B_1,B_2)=B_2$. We have 
$\p_{1!}\p_2^*=\p'_{1!}\p''_1{}^*\p'_{2!}\p''_2{}^*$. 
(Indeed, there are canonical maps 
$$\matrix
&& \hskip -4pt\fO_y &&\\
& \hskip -8pt\raise.6ex\hbox{${}^{\r_1}$}\hskip -8pt\swarrow &&\hskip
-4pt \searrow\hskip -5pt\raise.6ex\hbox{${}^{\r_2}$}& \\
\fO_s&&&& \hskip -8pt\fO_{y'}\\
\endmatrix$$
such that $\p_1=\p'_1\r_1$, $\p_2=\p'_2\r_2$ and 
$$\matrix
&&\hskip -8pt \fO_y &&\\
& \hskip -8pt\raise.6ex\hbox{${}^{\r_1}$}\hskip -8pt\swarrow &&\hskip
-4pt \searrow\hskip
-5pt\raise.6ex\hbox{${}^{\r_2}$}& \\
\fO_s&&&& \hskip -8pt\fO_{y'}\\
& \hskip -4pt \searrow\hskip -5pt\raise.6ex\hbox{${}^{\p''_1}$} &&
\hskip -8pt\raise.6ex\hbox{${}^{\p'_2}$}\hskip -6pt\swarrow &\\ \vspace{.5ex}
&&\hskip -8pt \cb&&
\endmatrix$$
is a cartesian diagram. Thus we have $\p_{1!}=\p'_{1!}\r_{1!}$,
$\p_2^*=\r_2^*\p'_{2!}$, $\p''_1{}^*\p'_{2!}=\r_{1!}\r_2^*$ hence
$\p_{1!}\p_2^*= \p'_{1!}\r_{1!}\r_2^*\p'_{2!}=
\p'_{1!}\p''_1{}^*\p'_{2!} \p''_2{}^*$.) Hence the
$(\p_{1!}\p_2^*\cs)^i$ are the end of a spectral sequence starting with 
$$(\p'_{1!}\p''_1{}^*(\p'_{2!}\p''_2{}^*\cs)^{h'})^h.\tag a$$
It is then enough to show that (a) with the isomorphism induced by $\Ps$ belongs to $\cc'$. By the induction 
hypothesis applied to $y'$ we see that $(\p'_{2!}\p''_2{}^*\cs)^{h'}$ with the isomorphism induced by $\Ps$ 
belongs to $\cc'$. We then use the fact that the lemma is already proved when $y=s$. This completes the proof.

\head 2. The category $\cc$\endhead
\subhead 2.1\endsubhead
In the remainder of this paper we will assume (as we may) that $\ph\colon G@>>>G$ in 1.1 satisfies in addition to the 
requirements 1.1(i)-(iv), the requirement that $\ph$ is the square of a Frobenius map $\ph_1\colon G@>>>G$ relative to
a split $\FF_{q'}$-rational structure on $G$ ($q'$ is a power of $p$ with $q'{}^2=q$) and 1.1(i)-(iv) are 
satisfied when $\ph$ is replaced by $\ph_1$.

From 1.6(a) applied to $\ph_1,\ph$ instead of $\ph,\ph^2$ we see that for any $\cl\in\fD$, there is a canonical 
choice of an isomorphism $\t^\cl\colon \ph^*\cl@>>>\cl$ such that for any $B\in\cb^{\ph_1}$, $\t^\cl_B\colon \cl_B@>>>\cl_B$ 
is the identity map. 

We set $\tph\colon =\s\ph=\ph\s\colon G@>>>G$. This is the Frobenius map for a (not necessarily split) $\FF_q$-rational 
structure on $G$. Now $\tph$ induces a  map $\cb@>>>\cb$ which is denoted again by $\tph$. Note that 
$\tph^2=\ph^2$.

Now $\fD^\s$ (see 0.1) is exactly the set of all $\cl\in\fD$ such that $\tph^*\cl$ is isomorphic to $\cl$ in 
$\cc_0$ (see 1.1(iii)). Let $\cl\in\fD^\s$ and let $\a^\cl$ be as in 0.1. We claim that 

{\it the two compositions 
$$\s^*\ph^*\cl@>\s^*(\t^\cl)>>\s^*\cl@>\a^\cl>>\cl\quad\text{ and
}\quad\ph^*\s^*\cl@>\ph^*(\a^\cl)>>\ph^*\cl@>\t^\cl>>\cl$$ 
coincide;}
\nl
we then denote these compositions by $\b^\cl\colon
\tph^*\cl@>>>\cl$. By the $K$-equivariance of $\cl$ it is enough  
to show that these compositions induce the same map on the stalks $(\tph^*\cl)_{\tph(B)}@>>>\cl_B$ for some 
$B\in[\cl]$. We take $B\in[\cl]$ such that $\ph_1(B)=B$. The two maps on the stalks are then the compositions
$$\cl_{\tph(B)}@>\t^\cl_{\s(B)}>>\cl_{\s(B)}@>\a^\cl_B>>\cl_B\qquad \text{
and }\qquad \cl_{\tph(B)}@>\a^\cl_{\ph(B)}>>\cl_{\ph(B)}@>\t^\cl_B>>\cl_B.$$
Since $\t^\cl_B=1$ and $\t^\cl_{\s(B)}=1$ (note 
that $\ph_1(\s(B))=\s(B)$ since $\ph_1\s=\s\ph_1$) both these compositions are equal to
$\cl_{\s(B)}@>\a^\cl_B>>\cl_B$. This proves our claim.

Let $\cc$ be the category whose objects are pairs $(\cs,\Ph)$ where
$\cs\in\cc_0$ and $\Ph$ is an isomorphism  
$\tph^*\cs@>>>\cs$ in $\cc_0$ such that, if $\Ps$ is the composition 
$\tph^{2*}\cs@>\tph^*(\Ph)>>\tph^*\cs@>\Ph>>\cs$, we have
$(\cs,\Ps)\in\cc'$ (note that  
$\tph^{2*}\cs=\ph^{2*}\cs$). A morphism between two objects
$(\cs,\Ph)$, $(\cs',\Ph')$ of $\cc$ is a morphism 
$\cs@>e>>\cs'$ in $\cc_0$ such that the diagram
$$\CD\tph^*\cs@>\Ph>>\cs\\
@V\tph^*e VV     @V e VV\\
\tph^*\cs'@>\Ph'>>\cs'\endCD$$
is commutative. 

Note that $(\cs,\Ph)\m(\cs,\Ph)^\Xi\colon =(\cs,\Ph\tph^*(\Ph))$ is a functor $\Xi\colon \cc@>>>\cc'$. Moreover, if 
$(\cs,\Ph)\in\cc$, then $(\cs,-\Ph)\in\cc$. 

(a) Let $\cl\in\fD^\s$. Then $(\cl,\b^\cl)\in\cc$. If $(\cl,\Ph)\in\cc$ then $\Ph=\pm\b^\cl$.
\nl
This is immediate.

(b) Now let $\cl\in\fD^\s$ and let $\b^{\cl\sh}\colon \b^*\cl^\sh@>>>\cl^\sh$ be the canonical extension of $\b^\cl$.
For any $h\in\NN$ let $\b^{\cl\sh2h}\colon \tph^*(\cl^{\sh 2h})@>>>\cl^{\sh 2h}$ be the isomorphism induced by 
$\b^{\cl\sh}$. We have $(\cl^{\sh 2h},\b^{\cl\sh2h})\in\cc$.
\nl
Indeed, $\b^{\cl\sh2h}\tph^*(\b^{\cl\sh2h})$ is equal to $\tt^{\cl\sh2h}$ in 1.6(b).

(c) Let $(\cs,\Ph)\in\cc$ and let $i\in\ZZ$. Let $\Ph^{(i)}\colon \tph^*(D^i\cs)@>>>D^i\cs$ be the inverse of the 
isomorphism $D^i\cs@>>>\tph^*(D^i\cs)$ induced by $D$. We have $(D^i\cs,\Ph^{(i)})\in\cc$.
\nl
Indeed, if $\Ps=\Ph\tph^*(\Ph)$ then $\Ps^{(i)}=\Ph^{(i)}\tph^*(\Ph^{(i)})$ is as in 1.6(c).
 
\subhead 2.2\endsubhead
For any $(\cs,\Ps)\in\cc'$ let $$\left(\sm0&\Ps\\1&0\esm\right)\colon
\tph^*(\cs\op\tph^*\cs)@>>>\cs\op\tph^*\cs$$  
(that is $\left(\sm0&\Ps\\1&0\esm\right)\colon
\tph^*\cs\op\tph^{2*}\cs@>>>\cs\op\tph^*\cs)$ be the isomorphism
sending $(a',b')$ to $(\Ps(b'),a')$. We set 
$$(\cs,\Ps)^\Th=(\cs\op\tph^*\cs,\left(\sm0&\Ps\\1&0\esm\right)).$$
Note that $(\cs,\Ps)^\Th\in\cc$. We have the following result.

(a) Let $(\cs,\Ph)\in\cc$. Then there exists an isomorphism in $\cc$
$$((\cs,\Ph)^\Xi)^\Th@>\si>>(\cs,\Ph)\op(\cs,-\Ph).$$ 
We define an isomorphism (in $\cc_0$) 
$$e\colon \cs\op\tph^*\cs@>\si>>\cs\op\cs, \qquad
(a,b)\m(a+\Ph(b),a-\Ph(b)).$$
Then 
$$\tph^*e\colon \tph^*\cs\op\tph^{2*}\cs@>\si>>\tph^*\cs\op\tph^*\cs, \qquad
(a',b')\m(a'+\tph^*(\Ph)(b'),a'-\tph^*(\Ph)(b')).$$
Define 
$$\Ph_0\colon \tph^*\cs\op\tph^{2*}\cs@>>>\cs\op\tph^*\cs, \qquad
(a',b')\m(\Ph\tph^*(\Ph)(b'),a').$$
Define 
$$\Ph_1\colon \tph^*\cs\op\tph^*\cs@>>>\cs\op\cs, \qquad
(c,d)\m(\Ph(c),-\Ph(d)).$$
We have 
$$\align e\Ph_0(a',b')&=e(\Ph\tph^*(\Ph)(b'),a')\\
&=(\Ph\tph^*(\Ph)(b')+\Ph(a'),\Ph\tph^*(\Ph)(b')-\Ph(a')),\endalign$$
$$\align\Ph_1\tph^*e(a',b')&=\Ph_1(a'+\tph^*(\Ph)(b'),a'-\tph^*\Ph(b'))\\&= 
(\Ph(a')+\Ph\tph^*(\Ph)(b'),-\Ph(a')+\Ph\tph^*(\Ph)(b')).\endalign$$
Thus $e\Ph_0=\Ph_1\tph^*e$ so that $e$ is an isomorphism as in (a); this proves (a).

\subhead 2.3\endsubhead
Let $\frak K(\cc)$ (resp. $\frak K(\cc')$) be the Grothendieck group of $\cc$ (resp. of $\cc'$). From the definitions we see
that the elements $(\cl,q^{2k}\tt^\cl)$ (with $\cl\in\fD$ and $k\in\ZZ$) form a $\ZZ$-basis of $\frak K(\cc')$. 
Moreover, the elements $(\cl,q^{2k}\tt^\cl)^\Th$ (with $\cl\in\fD-\fD^\s$ and $k\in\ZZ$) together with the 
elements $(\cl,\pm q^k\b^\cl)$ (with $\cl\in\fD^\s$ and $k\in\ZZ$) form a $\ZZ$-basis of $\frak K(\cc)$. Now $\Xi$ in 
2.1 defines a group homomorphism $\frak K(\cc)@>>>\frak K(\cc')$ denoted again by $\Xi$ and $\Th$ in 2.2 defines a group 
homomorphism $\frak K(\cc')@>>>\frak K(\cc)$ denoted again by $\Th$. From 2.2(a) we see that the image of this last 
homomorphism is the subgroup of $\frak K(\cc)$ spanned by the elements $(\cl,q^{2k}\tt^\cl)^\Th$ (with 
$\cl\in\fD-\fD^\s$ and $k\in\ZZ$) and by the elements $(\cl,q^k\b^\cl)+(\cl,-q^k\b^\cl)$ (with $\cl\in\fD^\s$ and
$k\in\ZZ$). It follows that 
$$M:=\frak K(\cc)/\Th(\frak K(\cc'))$$ 
has a $\ZZ$-basis consisting of the elements $(\cl,q^k\b^\cl)$ (with $\cl\in\fD^\s$ and $k\in\ZZ$). Note that in 
$M$ we have the equality $(\cl,-q^k\b^\cl)=-(\cl,q^k\b^\cl)$ (see 2.2(a)).

Let $\ca=\ZZ[u,u\i]$ where $u$ is an indeterminate. We regard $\frak K(\cc)$ as an $\ca$-module where 
$u^n(\cs,\Ph)=(\cs,q^n\Ph)$ for $n\in\ZZ$. We regard $\frak K(\cc')$ as an $\ca$-module where 
$u^n(\cs,\Ps)=(\cs,q^{2n}\Ph)$ for $n\in\ZZ$. Clearly, $\Xi\colon \frak K(\cc)@>>>\frak K(\cc')$ and 
$\Th\colon \frak K(\cc')@>>>\frak K(\cc)$ are
$\ca$-linear. Hence $M$ inherits from $\frak K(\cc)$ an $\ca$-module structure. Note that the elements 
$$a_\cl:=(\cl,\b^\cl)$$
(with $\cl\in\fD^\s$) form an $\ca$-basis of $M$.

\subhead 2.4\endsubhead
Clearly, there is a well-defined $\ZZ$-linear map $\DD\colon \frak K(\cc)@>>>\frak K(\cc)$ such that
$$\DD(\cs,\Ph)=\sum_{i\in\ZZ}(-1)^i(D^i\cs,\Ph^{(i)})$$ 
(notation of 2.1(c)) for any $(\cs,\Ph)\in\cc$. Similarly, there is a well-defined $\ZZ$-linear map 
$\DD\colon \frak K(\cc')@>>>\frak K(\cc')$ such that $$\DD(\cs,\Ps)=\sum_{i\in\ZZ}(-1)^i(D^i\cs,\Ps^{(i)})$$ 
(notation of 1.6(c)) for any $(\cs,\Ps)\in\cc'$. The homomorphism $\Th\colon \frak K(\cc')@>>>\frak K(\cc)$ is compatible with the 
maps $\DD$ hence $\DD\colon \frak K(\cc)@>>>\frak K(\cc)$ induces a $\ZZ$-linear map $M@>>>M$ denoted again by $\DD$. From the 
definitions we see that $\DD(u^n\x)=u^{-n}\DD(\x)$ for any $\x\in M$ and any $n\in\ZZ$.

\head 3. The category $\fC$\endhead
\subhead 3.1\endsubhead
In this section we specialize the definitions and results of \S1, \S2
in the context of 0.2. (In particular $G$ is 
replaced by $G\T G$.) In this case $\ph$ is replaced by
$\uph:=\ph\T\ph$ (hence $\tph$ is replaced by  
$\utph:=\ph\T\ph$). If $\fS\in\fC_0$ and $(B,B')\in\cb\T\cb$, let
$\fS_{B,B'}$ be the stalk of $\fS$ at $(B,B')$.  
For any $w\in W$ let $t^w\colon \uph^{2*}\SS_w@>>>\SS_w$ be the
isomorphism such that for any $(B,B')\in\fO_w$, the  
induced map on stalks $t^w_{B,B'}\colon
(\uph^{2*}\SS_w)_{B,B'}@>>>(\SS_w)_{B,B'}$ is the identity map
$\bbq@>>>\bbq$.  

Now $\cc'$ becomes $\fC'$; the objects of $\fC'$ are pairs $(\fS,\Ps)$
where $\fS\in\fC_0$ and $\Ps$ is an  
isomorphism $\uph^{2*}\fS@>\si>>\fS$ in $\fC_0$ such that for any
$(B,B')\in(\cb\T\cb)^{\uph}$, the eigenvalues  
of the induced map $\Ps_{B,B'}\colon
\fS_{B,B'}=\fS_{\ph^2B,\ph^2B'}@>>>\fS_{B,B'}$ are of the form $q^{2e}
(e\in\ZZ)$.  
For example, if $w\in W$ then $(\SS_w,t^w)\in\fC'$.

Moreover, $\cc$ 
becomes $\fC$; the objects of $\fC$ are pairs $(\fS,\Ph)$ where
$\fS\in\fC_0$ and $\Ph$ is an isomorphism  
$\utph^*\fS@>>>\fS$ in $\fC_0$ such that, if $\Ps$ is the composition
$$\utph^{2*}\fS@>\utph^*(\Ph)>>\utph^*\fS@>\Ph>>\fS,$$
we have $(\fS,\Ps)\in\fC'$ (note that 
$\utph^{2*}\fS=\uph^{2*}\fS$). Note that
$$(\fS,\Ph)\m(\fS,\Ph)^\Xi:=(\fS,\Ph\utph^*(\Ph))$$
is a functor $\Xi\colon \fC@>>>\fC'$. Moreover, if $(\fS,\Ph)\in\fC$, then
$(\fS,-\Ph)\in\fC$. Note that for $w\in W$ we have  
$\utph(\fO_w)=\fO_{\s(w)}$. For $w\in W^\s$ we denote by $\b^w\colon
\utph^*\SS_w@>>>\SS_w$ the isomorphism in $\fC_0$  
such that for any $(B,B')\in\fO_w$, the induced map on stalks 
$\b^w_{B,B'}\colon (\utph^*\SS_w)_{B,B'}@>>>(\SS_w)_{B,B'}$ is the
identity map $\bbq@>>>\bbq$. Clearly we have  
$(\SS_w,\b^w)\in\fC$.

Let $\b^{w\sh}\colon \utph^*\SS_w^\sh@>>>\SS_w^\sh$ be the isomorphism
which extends $\b^w$. For any $h\in2\NN$ let  
$\b^{w\sh2h}\colon \utph^*\SS_w^{\sh2h}@>>>\SS_w^{\sh2h}$ be the
isomorphism induced by $\b^{w\sh}$. We have  
$(\SS_w^{\sh2h},\b^{w\sh2h})\in\fC$. For any $(\fS,\Ps)\in\fC'$ let
$$\left(\sm0&\Ps\\1&0\esm\right)\colon
\utph^*(\fS\op\utph^*(\fS))@>>>\fS\op\utph^*(\fS)$$ 
(that is  
$$\left(\sm0&\Ps\\1&0\esm\right)\colon
\utph^*(\fS)\op\utph^{2*}(\fS)@>>>\fS\op\utph^*(\fS))$$ 
be the isomorphism given by $(a',b')\m(\Ps(b'),a')$. We set 
$$(\fS,\Ps)^\Th=(\fS\op\utph^*(\fS),\left(\sm0&\Ps\\1&0\esm\right)).$$
Note that $(\fS,\Ps)^\Th\in\fC$. For any $(\fS,\Ph)\in\fC$ there
exists an isomorphism in $\fC$
$$((\fS,\Ph)^\Xi)^\Th@>>>(\fS,\Ph)\op(\fS,-\Ph).$$ 

Let $\frak K(\fC)$ (resp. $\frak K(\fC')$) be the Grothendieck group of $\fC$ (resp. of $\fC'$); these are special cases of 
$\frak K(\cc),\frak K(\cc')$. Also, $\frak K(\fC)$, $\frak K(\fC')$ are naturally $\ca$-modules and $\Xi,\Th$ define $\ca$-linear maps
$\frak K(\fC)@>>>\frak K(\fC')$, $\frak K(\fC')@>>>\frak K(\fC)$ denoted again by $\Xi,\Th$. We set 
$\HH=\frak K(\fC)/\Th(\frak K(\fC'))$. This is 
naturally an $\ca$-module. It has an $\ca$-basis indexed by $W^\s$; to $w\in W^\s$ corresponds the element 
represented by $(\SS_w,\b^w)\in\fC$.

Let $\DD\colon \frak K(\fC)@>>>\frak K(\fC)$ be the $\ZZ$-linear map which is a special case of the map 
$\DD\colon \frak K(\cc)@>>>\frak K(\cc)$ in 
2.4. As in 2.4 this induces a $\ZZ$-linear endomorphism of $\HH=\frak K(\fC)/\Th(\frak K(\fC'))$ denoted again by $\DD$; it 
satisfies $\DD(u^n\x)=u^{-n}\DD(\x)$ for any $\x\in\HH$ and any $n\in\ZZ$.

\subhead 3.2\endsubhead
Let $y\in W$. We consider the diagram $\cb\T\cb@<\p_{13}<<\fO_y\T\cb@>\p_{23}>>\cb\T\cb$ where 
$\p_{13}(B,B',B'')=(B,B'')$, $\p_{23}(B,B',B'')=(B',B'')$. For $(\fS,\Ps)\in\fC'$ and $i\in\ZZ$ we have
$(\p_{13!}\p_{23}^*\fS)^i\in\fC_0$; we denote by 
$\Ps_i\colon \uph^{2*}(\p_{13!}\p_{23}^*\fS)^i@>>>(\p_{13!}\p_{23}^*\cs)^i$ the isomorphism induced by $\Ps$. 

\proclaim{Lemma 3.3} We have $((\p_{13!}\p_{23}^*\fS)^i,\Ps_i)\in\fC'$.
\endproclaim
The proof is similar to that of Lemma 1.8. We argue by induction on $l(y)$. If $y=1$ there is nothing to prove. 
Assume now that $l(y)=1$. By a standard argument we can assume that $(\fS,\Ps)=(\SS_w,t^w)$ (with $w\in W$). 
We consider the map $X@>c>>\cb\T\cb$ where 
$$X=\{(B,B',B'')\in\cb\T\cb\T\cb;(B,B')\in\fO_y,(B',B'')\in\fO_w\},$$
$c(B,B',B'')=(B,B'')$. For $i\in\ZZ$ we have an obvious isomorphism 
$$\t_i\colon \uph^{*2}(c_!\bbq)^i@>>>(c_!\bbq)^i.$$
Clearly, $(c_!\bbq)^i\in\fC_0$. It is enough to show that $((c_!\bbq)^i,\t_i)\in\fC'$. This is easily verified
since any fibre of $c$ is either a point or $\kk$ or $\kk^*$.

Next we assume that $l(y)\ge2$ and that the result is known for elements of length $<\ul(y)$. We can find 
$s\in S$ and $y'\in W$ such that $y=sy'$, $\ul(y)=\ul(y')+1$. We have a diagram
$$\cb\T\cb@<\p'_{13}<<\fO_s\T\cb@>\p'_{23}>>\cb\T\cb@<\p''_{13}<<\fO_{y'}\T\cb@>\p''_{23}>>\cb\T\cb$$ 
where $\p'_{13},\p''_{13}$ are given by $(B,B',B'')\m(B,B'')$ and $\p'_{23},\p''_{23}$ are given by 
$(B,B',B'')\m(B',B'')$. We have $\p_{13!}\p_{23}^*=\p'_{13!}\p'_{23}{}^*\p''_{13!}\p''_{23}{}^*$. Hence the \lb
$(\p_{13!}\p_{23}^*\fS)^i$ are the end of a spectral sequence starting with
$$(\p'_{13!}\p'_{23}{}^*(\p''_{13!}\p''_{23}{}^*\fS)^{h'})^h.\tag a$$
It is then enough to show that (a) with the isomorphism induced by $\Ps$ belongs to $\fC'$. By the induction 
hypothesis applied to $y'$ we see that $(\p''_{13}\p''_{23}{}^*\fS)^{h'}$ with the isomorphism induced by $\Ps$ 
belongs to $\fC'$. We then use the fact that the lemma is already proved when $y=s$. This completes the proof.

\head 4. $\HH$ as an algebra and $M$ as an $\HH$-module\endhead
\subhead 4.1 \endsubhead
We return to the setup in 2.1. We consider the diagram $\cb@<\p_1<<\cb\T\cb@>\p_2>>\cb$ where $\p_1,\p_2$ are the
first and second projection. For any $\fS\in\fC_0,\cs\in\cc_0$ and any $i\in\ZZ$ we set
$$\fS\odot_i\cs=(\p_{1!}(\fS\ot(\p_2^*\cs)))^i\in\cc_0.$$
If $t\colon \fS@>>>\fS'$ and $t'\colon \cs@>>>\cs'$ are isomorphisms in $\fC_0$ and $\cc_0$ then the induced isomorphism 

(a) $\fS\odot_i\cs@>>>\fS'\odot_i\cs'$ (in $\cc_0$) is denoted by $t\odot_it'$. 
\nl
We show:

(b) If $(\fS,\Ps)\in\fC'$ and $(\cs,\Ps')\in\cc'$ then $(\fS\odot_i\cs,\Ps\odot_i\Ps')\in\cc'$.
\nl
By a standard argument we may reduce the general case to the case where for some $w\in W$, $\fS=\SS_w$, $\Ps=t_w$.
In this case the result follows from 1.8.

We show:

(c) If $(\fS,\Ph)\in\fC$ and $(\cs,\Ph')\in\cc$ then $(\fS\odot_i\cs,\Ph\odot_i\Ph')\in\cc$.
\nl
It is enough to show that $(\fS\odot_i\cs,(\Ph\odot_i\Ph')\tph^*(\Ph\odot_i\Ph'))\in\cc'$ or that \lb
$(\fS\odot_i\cs,(\Ph\utph^*\Ph)\odot_i(\Ph'\tph^*\Ph'))\in\cc'$. This follows from (b).

In the setup of (c) it makes sense to define
$$(\fS,\Ph)*(\cs,\Ph')=\sum_i(-1)^i(\fS\odot_i\cs,\Ph\odot_i\Ph')\in\frak K(\cc).$$
This gives rise to an $\ca$-bilinear pairing
$$\frak K(\fC)\T\frak K(\cc)@>>>\frak K(\cc).\tag d$$
Let $(\fS,\Ps)\in\fC'$ and $(\cs,\Ph')\in\cc$. Then $(\fS,\Ps)^\Th\in\fC$. We have the following result.
$$(\fS,\Ps)^\Th*(\cs,\Ph')=\sum_i(-1)^i(\fS\odot_i\cs,\Ps\odot_i(\Ph'\tph^*(\Ph')))^\Th.\tag e$$
It is enough to show that for any $i\in\ZZ$ the following diagram is commutative
$$\CD\utph^*\fS\odot_i\tph^*\cs\op\utph^{2*}\fS\odot_i\tph^{2*}\cs@>A>>
\fS\odot_i\cs\op\utph^*\fS\odot_i\tph^*\cs\\
@VCVV                      @VC'VV\\
\utph^*\fS\odot_i\tph^*\cs\op\utph^{2*}\fS\odot_i\tph^*\cs@>A'>>\fS\odot_i\cs\op\utph^*\fS\odot_i\cs\endCD$$
where 
$$\align &A=\left(\sm0&\Ps\odot_i\Ph'\tph^*(\Ph')\\1&0\esm\right),\qquad
A'=\left(\sm0&\Ps\odot_i\Ph'\\1\odot_i\Ph'&0\esm\right),\\&
C=\left(\sm1\odot_i1&0\\0&1\odot_i\tph^*(\Ph')\esm\right),\qquad
C'=\left(\sm1\odot_i1&0\\0&1\odot_i\Ph'\esm\right).\endalign$$
This is immediate.

Now let $(\fS,\Ph)\in\fC$ and $(\cs,\Ps')\in\cc'$. Then $(\cs,\Ps')^\Th\in\cc$. A proof analogous to that of (e) 
shows that
$$(\fS,\Ph)*(\cs,\Ps')^\Th=\sum_i(-1)^i(\fS\odot_i\cs,\Ph\utph^*(\Ph)\odot_i\Ps'))^\Th.\tag f$$
From (e),(f) we see that the pairing (d) factors through an $\ca$-bilinear pairing 
$$\HH\T M@>>>M.\tag g$$

\subhead 4.2\endsubhead 
Let $\p_{ab}\colon \cb\T\cb\T\cb@>>>\cb\T\cb$ ($a,b$ is $12$ or $23$
or $13$) be the projection to the $a,b$ factors. 
For any $\fS,\fS'\in\fC_0$ and any $i\in\ZZ$ we set
$$\fS\odot_i\fS'=\p_{13!}(\p_{12}^*\fS\ot\p_{23}^*\fS')\in\fC_0.$$
If $t\colon \fS@>>>\fS_1$ and $t'\colon \fS'@>>>\fS'_1$ are isomorphisms in $\fC_0$ then the induced isomorphism 

(a) $\fS\odot_i\fS'@>>>\fS_1\odot_i\fS'_1$ (in $\fC_0$) is denoted by
$t\odot_it'$.  
\nl
We show:

(b) If $(\fS,\Ps)\in\fC'$ and $(\fS',\Ps')\in\fC'$ then $(\fS\odot_i\fS',\Ps\odot_i\Ps')\in\fC'$.
\nl
By a standard argument we may reduce the general case to the case where $\fS=\SS_w$ for some $w\in W$ and 
$\Ps=t^w$ is as in 3.1. In this case the result follows from 3.3.

We show:

(c) If $(\fS,\Ph)\in\fC$ and $(\fS',\Ph')\in\fC$ then
$(\fS\odot_i\fS',\Ph\odot_i\Ph')\in\fC$. 
\nl
It is enough to show that 
$$(\fS\odot_i\fS',(\Ph\odot_i\Ph')\utph^*(\Ph\odot_i\Ph'))\in\fC'$$
or that 
$$(\fS\odot_i\fS',(\Ph\utph^*(\Ph))\odot_i(\Ph'\utph^*(\Ph')))\in\fC'.$$
This follows from (b).

In the setup of (c) it makes sense to define
$$(\fS,\Ph)*(\fS',\Ph')=\sum_i(-1)^i(\fS\odot_i\fS',\Ph\odot_i\Ph')\in\frak
K(\fC).$$ 
This gives rise to a $\ZZ$-bilinear pairing
$$\frak K(\fC)\T\frak K(\fC)@>>>\frak K(\fC).\tag d$$
A standard argument shows that the pairing (d) defines an associative
$\ca$-algebra structure on $\frak K(\fC)$ with 
unit element represented by $(\SS_1,\b^1)$ and that 

(e) {\it the pairing 4.1(d) defines a (unital) $\frak K(\fC)$-module
structure on the $\ca$-module $\frak K(\cc)$.} 

Let $(\fS,\Ps)\in\fC'$ and $(\fS',\Ph')\in\fC$. Then
$(\fS,\Ps)^\Th\in\fC$. We have the following analogue of 4.1(e): 
$$(\fS,\Ps)^\Th*(\fS',\Ph')=\sum_i(-1)^i(\fS\odot_i\fS',\Ps\odot_i(\Ph'\utph^*(\Ph')))^\Th.\tag f$$
It is enough to show that for any $i\in\ZZ$ the following diagram is commutative
$$\CD\utph^*\fS\odot_i\utph^*\fS'\op\utph^{2*}\fS\odot_i\utph^{2*}\fS'@>A>>
\fS\odot_i\fS'\op\utph^*\fS\odot_i\utph^*\fS'\\
@VCVV                      @VC'VV\\
\utph^*\fS\odot_i\tph^*\fS'\op\utph^{2*}\fS\odot_i\utph^*\fS'@>A'>>\fS\odot_i\fS'\op\utph^*\fS\odot_i\fS'\endCD$$
where  
$$\align& A=\left(\sm0&\Ps\odot_i\Ph'\utph^*(\Ph')\\1&0\esm\right),\qquad
A'=\left(\sm0&\Ps\odot_i\Ph'\\1\odot_i\Ph'&0\esm\right),\\&
C=\left(\sm1\odot_i1&0\\0&1\odot_i\utph^*(\Ph')\esm\right),\qquad
C'=\left(\sm1\odot_i1&0\\0&1\odot_i\Ph'\esm\right).\endalign$$
This is immediate.

Now let $(\fS,\Ph)\in\fC$ and $(\fS',\Ps')\in\fC'$. Then $(\fS',\Ps')^\Th\in\fC$. A proof analogous to that of 
(f) shows that
$$(\fS,\Ph)*(\fS',\Ps')^\Th=\sum_i(-1)^i(\fS\odot_i\fS',\Ph\utph^*(\Ph)\odot_i\Ps')^\Th.\tag g$$
From (f), (g) we see that $\Th(\frak K(\fC'))$ is a two-sided ideal of
$\frak K(\fC)$ hence $\HH$ inherits from  $\frak K(\fC)$ a 
structure of associative $\ca$-algebra with $1$. Hence the pairing $\HH\T M@>>>M$ in 4.1(g) makes the $\ca$-module $M$ into a (unital) $\HH$-module (see (e)).

\subhead 4.3\endsubhead
The subgroup $W^\s$ of $W$ is itself a Coxeter group with standard
generators $w_\o$ where $\o$ runs over the set  
$\bS$ of $\s$-orbits in $S$ and $w_\o$ is the longest element in the
subgroup of $W$ generated by the elements in 
$\o$. 

\subhead 4.4\endsubhead
In the setup of 3.1, let $s$ be an odd integer greater than or equal
to $1$. Note that $\tph^s$ is the Frobenius map for an 
$\FF_{q^s}$-rational structure on $G$ and 
$$\utph^s\colon \cb\T\cb@>>>\cb\T\cb, \qquad
(B,B')\m(\tph^s(B),\tph^s(B'))$$ 
is the Frobenius map for an $\FF_{q^s}$-rational structure on
$\cb\T\cb$. Now $G^{\tph^s}$ acts on 
$(\cb\T\cb)^{\utph^s}$ (by restriction of the $G$-action on $\cb\T\cb$) and from Lang's theorem we see that the 
$G^{\tph^s}$-orbits on $(\cb\T\cb)^{\utph^s}$ are exactly the sets $\fO_w^{\utph^s}$ where $w$ runs through
the set of fixed points of $\s^s$ on $W$ (which is the same as $W^\s$ since $s$ is odd).

Let $\cg_s$ be the $\bbq$-vector space consisting of all functions $f\colon (\cb\T\cb)^{\utph^s}@>>>\bbq$ which are 
constant on each $G^{\tph^s}$-orbit on $(\cb\T\cb)^{\utph^s}$. For any $w\in W^\s$ let $f_{w,s}$ be the function 
on $(\cb\T\cb)^{\utph^s}$ which is equal to $1$ on $\fO_w^{\utph^s}$ and is equal to $0$ on $\fO_{w'}^{\utph^s}$ 
for $w'\in W^\s-\{w\}$. Note that $\{f_{w,s};w\in W^\s\}$ is a $\bbq$-basis of $\cg_s$.

For $f,f'\in\cg_s$ we define $f*f'\in\cg_s$ by
$$(f*f')(B,B'')=\sum_{B'\in\cb^{\tph^s}}f(B,B')f'(B',B'').$$
This defines on $\cg_s$ a structure of associative $\bbq$-algebra with unit element $f_{1,s}$. From 
\cite{\IW}, \cite{\MA} it follows that the following relations hold in this algebra:
$$f_{w,s}*f_{w',s}=f_{ww',s}\text{ if }w,w'\in W^\s, \ul(ww')=\ul(w)+\ul(w');\tag a$$
$$f_{w_\o,s}*f_{w_\o,s}=q^{s\ul(w_\o)}f_{1,s}+(q^{s\ul(w_\o)}-1)f_{w_\o,s}\text{ if }\o\in\bS\tag b$$ 
(see 4.3).

\subhead 4.5\endsubhead
In the setup of 4.4, for any $(\fS,\Ph)\in\fC$ we define an isomorphism $\Ph_s\colon \utph^{s*}\fS@>>>\fS$ as the 
composition
$$\utph^{s*}\fS@>\utph^{(s-1)*}(\Ph)>>\utph^{(s-1)*}\fS@>>>\do@>>>\utph^*\fS@>\Ph>>\fS.$$
For any $(B,B')\in(\cb\T\cb)^{\utph^s}$, $\Ph_s$ induces a linear isomorphism $\fS_{B,B'}@>>>\fS_{B,B'}$ whose 
trace is denoted by $\c_{s;\fS,\Ph}(B,B')$. For any $g\in G^{\tph^s}$ we have 
$$\c_{s;\fS,\Ph}(gBg\i,gB'g\i)=\c_{s;\fS,\Ph}(B,B').$$
Hence the function $(B,B')\m\c_{s,\fS,\Ph}(B,B')$ belongs to $\cg_s$. Note that $(\fS,\Ph)\m\c_{s,\fS,\Ph}$ 
defines a group homomorphism $\frak K(\fC)@>>>\cg_s$. From the definitions we see that the kernel of this homomorphism 
contains $\Th(\frak K(\fC'))$ hence we get an induced group homomorphism $\vt_s\colon \HH@>>>\cg_s$ such that 
$\vt_s(u^n\x)=q^{ns}\vt_s(\x)$ for all $\x\in\HH$, $n\in\ZZ$ and such that $\vt_s(\SS_w,\b^w)=f_{w,s}$ for any 
$w\in W^\s$. It follows that $\vt_s$ induces an isomorphism of $\bbq$-vector spaces
$$\bar\vt_s\colon \bbq\ot_\ca\HH@>>>\cg_s$$
where $\bbq$ is regarded as an $\ca$-algebra by $u\m q^s$.

Now let $(\fS,\Ph)\in\fC$, $(\fS',\Ph')\in\fC$. Using the definitions and Grothendieck's ``faisceaux-fonctions" 
dictionary we see that
$$\vt_s(\fS,\Ph)*\vt_s(\fS',\Ph')=\vt_s((\fS,\Ph)*(\fS',\Ph')).$$
In other words, 

(a) $\vt_s\colon \HH@>>>\cg_s$ is a ring homomorphism (hence $\bar\vt_s\colon \bbq\ot_\ca\HH@>>>\cg_s$ is an algebra
isomorphism).

\proclaim{Lemma 4.6} The homomorphism $\vt\colon \HH@>>>\op_{s\in\{1,3,5,\do\}}\cg_s$ (with components $\vt_s$) is 
injective.
\endproclaim
Let $\x\in\HH$ be such that $\vt(\x)=0$. We can write uniquely 
$$\x=\sum_{w\in W^\s}c_w(u)(\SS_w,\b^w)$$
where $c_w(u)\in\ca$. Using our assumption we deduce that $\sum_{w\in W^\s}c_w(q^s)f_{w,s}=0$ for any 
$s\in\{1,3,\do\}$. Since $f_{w,s}$ are linearly independent in $\cg_s$
we deduce that $c_w(q^s)=0$ for any  
$s\in\{1,3,\do\}$ and any $w\in W^\s$. Hence $c_w(u)=0$ for any $w\in W^\s$ so that $\x=0$. The lemma is proved.

\subhead 4.7\endsubhead
For any $w\in W^\s$ let $T_w=(\SS_w,\b^w)\in\HH$. The following
identities hold in $\HH$: 
$$T_wT_{w'}=T_{ww'}\text{ if }w,w'\in W^\s,
\ul(ww')=\ul(w)+\ul(w');\tag a$$ 
$$T_{w_\o}^2=u^{\ul(w_\o)}T_1+(u^{\ul(w_\o)}-1)T_{w_\o}\text{ if
}\o\in\bS\tag b$$  
(see 4.3). Indeed, using Lemma 4.6 it is enough to show that these
identities hold after applying $\vt_s$ for any 
$s\in\{1,3,\do\}$. But this follows from 4.4(a), (b). (Note that (a)
can be proved also directly from the definitions.)

This shows that $\HH$ with its $\ca$-basis $\{T_w;w\in W^\s\}$ is a
(quasisplit) Iwahori-Hecke algebra. 

\subhead 4.8\endsubhead 
Let $\o\in\bS$. Let $(\fS_{\le w_\o},\Ph)\in\fC$ be such that  
$$\fS_{\le w_\o}|_{\ov{\fO_{w_\o}}}=\bbq,\qua \fS_{\le
w_\o}|_{(\cb\T\cb)-\ov{\fO_{w_\o}}}=0$$ 
and for any $(B,B')\in\ov{\fO_{w_\o}}$, $\Ph$ induces the identity map
from  $(\fS_{\le w_\o})_{\tph(B),\tph(B')}=\bbq$ to $(\fS_{\le
w_\o})_{B,B'}=\bbq$. From the definitions we see that  
$(\fS_{\le w_\o},\Ph)$ represents the element $T_{w_\o}+T_1$ of
$\HH$. (The only elements of $W^\s$ which are  
contained in the subgroup of $W$ generated by the elements in $w_\o$
are $w_\o$ and $1$.)  

Now consider the diagram
$\cb\T\cb@<\p_{13}<<\ov{\fO_{w_\o}}\T\cb@>\p_{23}>>\cb\T\cb$ where  
$\p_{13},\p_{23}$ are the restrictions of the maps with the same name in 4.2.
For any $(\fS',\Ph')\in\fC$ and any $i\in\ZZ$ let
$\fS'(i)=(\p_{13!}\p_{23}^*\fS')^i$ and let $\Ph'(i)\colon \utph^*\fS'(i)@>>>\fS'(i)$ be the isomorphism induced by
$\Ph'$. We have $(\fS'(i),\Ph'(i))\in\fC$ (a special case of 4.2(c)).
We set $\th_\o(\fS',\Ph')=\sum_i(-1)^i(\fS'(i),\Ph'(i))$. We can view $\th_\o$ as an $\ca$-linear map
$\frak K(\fC)@>>>\frak K(\fC)$. From the definitions we have
$\th_\o(\x)=(\fS_{\le w_\o},\Ph)\x$ in the algebra $\frak K(\fC)$.
From the known properties of Verdier duality we have that
$$\DD(\th_\o(\x))=u^{-\ul(w_\o)}\th_\o(\DD(\x))\text{ for all }\x\in\frak K(\fC).   \tag a$$
(We use that $\p_{13}$ is proper and $\p_{23}$ is smooth with connected fibres of dimension $\ul(w_\o)$.) 
Clearly $\th_\o$ induces an $\ca$-linear map $\HH@>>>\HH$ denoted again by $\th_\o$ such that
in the algebra $\HH$ we have 
$$\th_\o(\x)=(T_{w_\o}+T_1)\x.$$
Applying $\DD$ to both sides and using (a) we obtain for any $\x\in\HH$:
$$u^{-\ul(w_\o)}(T_{w_\o}+T_1)\DD(\x)=\DD((T_{w_\o}+T_1)\x).\tag b$$
Now let $\,\bar{}\colon \HH@>>>\HH$ be the unique ring homomorphism such that $\ov{u^nT_w}=u^{-n}T_{w\i}\i$ for any 
$w\in W^\s, n\in\ZZ$. Note that $\ov{T_{w_\o}+T_1}=u^{-\ul(w_\o)}(T_{w_\o}+T_1)$. Hence (b) implies
$\DD(h\x)=\ov{h}\DD(\x)$ whenever $h=T_{w_\o}+T_1$, $\x\in\HH$. We have also $\DD(u^n\x)=\ov{u^n}\DD(\x)$ for
$n\in\ZZ$, $\x\in\HH$. Since the elements $h$ as above and $u^n$ ($n\in\ZZ$) generate the ring $\HH$ it follows 
that $\DD(h)=\ov{h}\DD(T_1)$ for any $h\in\HH$. From the definitions we have $\DD(T_1)=u^{-\nu}T_1$ where 
$\nu=\dim\cb$. It follows that
$$\DD(h)=u^{-\nu}\ov{h}\text{ for any }h\in\HH.\tag c$$

Now consider the diagram $\cb@<\p_1<<\ov{\fO_{w_\o}}@>\p_2>>\cb$ where
$\p_1,\p_2$ are the first and second  
projection. For any $(\cs,\Ph')\in\cc$ and any $i\in\ZZ$ let
$\cs(i)=(\p_{1!}\p_2^*\cs)^i$ and let  
$\Ph'(i)\colon \tph^*\cs(i)@>>>\cs(i)$ be the isomorphism induced by
$\Ph'$. We have $(\cs(i),\Ph'(i))\in\cc$ (a special case of 4.1(c)). We set
$\th_\o(\cs,\Ph')=\sum_i(-1)^i(\cs(i),\Ph'(i))$. We can view $\th_\o$
as an  
$\ca$-linear map $\frak K(\cc)@>>>\frak K(\cc)$. From the definitions
we have $\th_\o(\x)=(\fS_{\le w_\o},\Ph)\x$ in the  
$\frak K(\fC)$-module $\frak K(\cc)$. From the known properties of
Verdier duality we have that 
$$\DD(\th_\o(\x))=u^{-\ul(w_\o)}\th_\o(\DD(\x))\text{ for all
}\x\in\frak K(\cc).\tag d$$ 
(We use that $\p_1$ is proper and $\p_2$ is smooth with connected fibres of dimension $\ul(w_\o)$.) Clearly,
$\th_\o$ induces an $\ca$-linear map $M@>>>M$ denoted again by $\th_\o$ such that in the $\HH$-module $M$ we have 
$$\th_\o(\x)=(T_{w_\o}+T_1)(\x).$$
Applying $\DD$ to both sides and using (d) we obtain
$$u^{-\ul(w_\o)}(T_{w_\o}+T_1)(\DD(\x))=\DD((T_{w_\o}+T_1)(\x)).\tag e$$
Thus, $\DD(h\x)=\ov{h}\DD(\x)$ if $h=T_{w_\o}+T_1$, $\x\in M$. We have also $\DD(u^n\x)=u^{-n}\DD(\x)$ for 
$n\in\ZZ$, $\x\in M$. Since the elements $h$ as above and $u^n$ ($n\in\ZZ$) generate the ring $\HH$, it follows 
that 
$$\DD(h\x)=\bar h\DD(\x)\text{ for any }h\in\HH,\x\in M.\tag f$$

\subhead 4.9\endsubhead 
Let $\o\in\bS$ and let $\cl\in\fD^\s$. In the $\HH$-module $M$ we have
$$T_{w_\o}(\cl,\b^\cl)=\sum_{\cl'\in\fD^\s}f_{\o;\cl',\cl}(\cl',\b^{\cl'})$$
where $f_{\o;\cl',\cl}\in\ca$ are well defined. We want to make the
quantities $f_{\o;\cl',\cl}$ as explicit as  
possible. Let $\co=[\cl]$. Let $X=\{(B,B')\in\fO_{w_\o};B'\in\co\}$;
define 
$$\et\colon X@>>>X, \qquad (B,B')\m(\tph(B),\tph(B'))$$
(the Frobenius map for an $\FF_q$-rational structure on $X$). The
local system  $\ucl$ on $\co$ pulls back under the second projection to a local
system $\tucl$ on $X$ and   
$\b^\cl\colon \tph^*\cl@>>>\cl$ induces an isomorphism
$\ti\b^\cl\colon \et^*\tucl@>>>\tucl$. Let $\p\colon X@>>>\cb$ be the
first projection. For any $i\in\NN$ we have $(\p_!\tucl)^i\in\cc_0$
and $\ti\b^\cl$ induces an isomorphism 
$\tph^*(\p_!\tucl)^i@>>>(\p_!\tucl)^i$ (in $\cc_0$) denoted by $\b^{(i)}$. Note that $\p(X)$ is a union of 
$K$-orbits.  For any $K$-orbit $\co'$ contained in $\p(X)$ we have
canonically 
$$(\p_!\tucl)^i= \bigoplus_{\cl'\in\fD^\s;[\cl']=\co'}
V_{\cl'}((\p_!\tucl)^i) \ot\ucl'$$  
where $V_{\cl'}((\p_!\tucl)^i)$ are finite-dimensional $\bbq$-vector
spaces. Moreover in terms of this decomposition we have
$\b^{(i)}=\op_{\cl'\in\fD^\s;[\cl']=\co'}c_{\cl',i}\ot\b^{\cl'}$ where
$c_{\cl',i}$ is an automorphism of $V_{\cl'}((\p_!\tucl)^i)$ with all
eigenvalues of the form $\pm q^e$, $(e\in\ZZ)$; the dimension of the  
$(\pm q^e)$-eigenspace is denoted by $c_{\cl',i;\pm q^e}$. For any $\cl'\in\fD^\s$ such that 
$[\cl']\not\sub\p(X)$ we have $f_{\o;\cl',\cl}=0$. From the definitions, for any $\cl'\in\fD^\s$ such that 
$[\cl']\sub\p(X)$ we have
$$f_{\o;\cl',\cl}=\sum_{i\in\NN,e\in\ZZ}(-1)^i(c_{\cl',i;q^e}-c_{\cl',i;-q^e})u^e.$$
Note that the polynomials $\sum_{i\in\NN,e\in\ZZ}(-1)^i
(c_{\cl',i;q^e}+c_{\cl',i;-q^e})u^e$ are matrix  
coefficients of the action of $T_{w_\o}$ in a module \cite{\LV, 1.7}
over the split Hecke algebra associated to $W$; they can in principle
be calculated by iteration from \cite{\LV, Lemma 3.5}. 

\head 5. The elements $\fA_\cl\in M$\endhead 
\subhead 5.1\endsubhead 
For $\cl,\cl'\in\fD^\s$ we say that $\cl'\preceq\cl$ if

(i) $[\cl']\sub\ov{[\cl]}$ and

(ii) if $[\cl']=[\cl]$ then $\cl'=\cl$.
\nl
This defines a partial order on $\fD^\s$. We write $\cl'\prec\cl$ if $\cl'\preceq\cl,\cl'\ne\cl$. We have the 
following result.

\proclaim{Theorem 5.2} Let $\cl\in\fD^\s$. The polynomials $P^\s_{\cl',\cl}$ ($\cl'\in\fD^\s$) in 0.1 are
characterized by the following properties:
$$\sum_{\cl'\in\fD^\s}P^\s_{\cl',\cl}(u\i)\DD(a_{\cl'})=
u^{-\dim[\cl]}\sum_{\cl'\in\fD^\s}P^\s_{\cl',\cl}(u)a_{\cl'}\text{ in }M;\tag a$$
$$P^\s_{\cl',\cl}=0\text{ if }\cl'\not\preceq\cl;\tag b$$
$$\deg P^\s_{\cl',\cl}\le(\dim[\cl]-\dim[\cl']-1)/2\text{ if }\cl'<\cl\text{ and }P^\s_{\cl,\cl}=1.\tag c$$
\endproclaim
Let $\t^{\cl\sh}\colon \ph^*\cl^\sh@>>>\cl^\sh$ be the canonical extension of $\t^\cl$. Replacing $\s,\a^\cl,\a^{\cl'}$
by $\ph,\t^\cl,\b^{\cl'}$ in the definition of 
$\a^{\cl\sh2h;\cl'}\colon V_{\cl'}(\cl^{\sh2h})@>>>V_{\cl'}(\cl^{\sh2h})$ in 0.1 (here $h\in\NN$) we obtain a vector 
space isomorphism 
$$\t^{\cl\sh2h;\cl'}\colon V_{\cl'}(\cl^{\sh2h})@>>>V_{\cl'}(\cl^{\sh2h}).$$
From 1.6(b) applied to $\ph_1,\ph$ instead of $\ph,\ph^2$ we see that 

(d) {\it for any $h\in\NN$, $\t^{\cl\sh2h;\cl'}$ is equal to $q^h$
times a unipotent linear map.} 
\nl
Replacing $\s,\a^\cl,\a^{\cl'}$ by $\tph,\b^\cl,\b^{\cl'}$ in the definition of
$$\a^{\cl\sh2h;\cl'}\colon
V_{\cl'}(\cl^{\sh2h})@>>>V_{\cl'}(\cl^{\sh2h})$$ 
in 0.1 (here $h\in\NN$) we obtain a vector  
space isomorphism  
$$\b^{\cl\sh2h;\cl'}\colon V_{\cl'}(\cl^{\sh2h})@>>>V_{\cl'}(\cl^{\sh2h}).$$
Let $h\in\NN$. From the definitions we have 
$$\b^{\cl\sh2h;\cl'}=\t^{\cl\sh2h;\cl'}\a^{\cl\sh2h;\cl'}=\a^{\cl\sh2h;\cl'}\t^{\cl\sh2h;\cl'}.$$
For $\e=\pm1$ let $V_{\cl',\e}(\cl^{\sh2h})$ be the $\e$-eigenspace of  $\a^{\cl\sh2h;\cl'}$ (which has square 
$1$). We deduce that
$$\align&(V_{\cl'}(\cl^{\sh2h})\ot\cl',\b^{\cl\sh2h;\cl'}\ot\b^{\cl'})\\&=
(V_{\cl',1}(\cl^{\sh2h})\ot\cl',t^{\cl\sh2h;\cl'}\ot\b^{\cl'})+
(V_{\cl',-1}(\cl^{\sh2h})\ot\cl',-t^{\cl\sh2h;\cl'}\ot\b^{\cl'})\endalign$$
in $\frak K(\cc)$. Using (d) we deduce
$$\align&(V_{\cl'}(\cl^{\sh2h})\ot\cl',\b^{\cl\sh2h;\cl'}\ot\b^{\cl'})\\&=
\dim V_{\cl',1}(\cl^{\sh2h})u^h(\cl',\b^{\cl'})+\dim V_{\cl',-1}(\cl^{\sh2h})u^h(\cl',-\b^{\cl'})\endalign$$
in $\frak K(\cc)$ hence
$$\align&(V_{\cl'}(\cl^{\sh2h})\ot\cl',\b^{\cl\sh2h;\cl'}\ot\b^{\cl'})\\&=
u^h\tr(\a^{\cl\sh2h;\cl'}\colon V_{\cl'}(\cl^{\sh2h})@>>>V_{\cl'}(\cl^{\sh2h}))(\cl',\b^{\cl'})
=u^hP^\s_{\cl',\cl;h}a_{\cl'}\endalign$$
in $M$ (notation of 2.3). Using the definitions we have
$$(\cl^{\sh2h},\b^{\cl\sh2h})
=\sum_{\cl'\in\fD^\s;\cl'\preceq\cl}(V_{\cl'}(\cl^{\sh 2h})\ot\cl',\b^{\cl\sh2h;\cl'}\ot\b^{\cl'})$$
in $M$. Hence 
$$(\cl^{\sh2h},\b^{\cl\sh2h})=\sum_{\cl'\in\fD^\s;\cl'\preceq\cl}u^hP^\s_{\cl',\cl;h}a_{\cl'}$$
in $M$. Thus, setting
$$\fA_\cl=\sum_{h\in\NN}(\cl^{\sh2h},\b^{\cl\sh2h})\in M$$
we have 
$$\fA_\cl=\sum_{\cl'\in\fD^\s;\cl'\preceq\cl}P^\s_{\cl',\cl}a_{\cl'}\text{ in }M.$$
From the definition of $\cl^\sh$ we see that (b), (c) hold and that
$$\DD(\fA_\cl)=u^{-\dim[\cl]}\fA_\cl\text{ in }M.$$ 
Thus (a) holds. From (a), (b), (c) we deduce that for any
$\cl'\preceq\cl$ we have 
$$\DD(a_{\cl'})=\sum_{\cl''\in\fD^\s;\cl''\preceq\cl'}
\r_{\cl'',\cl'}a_{\cl''} \tag d$$  
where $\r_{\cl'',\cl'}\in\ca$ and $\r_{\cl',\cl'}=u^{-\dim[\cl']}$.

To complete the proof it is enough to show that if
$x_{\cl''}\in\ZZ[u]$ are defined for $\cl''\preceq\cl$ and  
satisfy $x_\cl=0$, $\deg x_{\cl''}\le(\dim[\cl]-\dim[\cl'']-1)/2$ if $\cl''<\cl$ and
$$\sum_{\cl'\in\fD^\s;\cl'\preceq\cl}x_{\cl'}(u\i)\sum_{\cl''\in\fD^\s;\cl''\preceq\cl'}\r_{\cl'',\cl'}a_{\cl''}
=u^{-\dim[\cl]}\sum_{\cl''\in\fD^\s;\cl''\preceq\cl}x_{\cl''}(u)a_{\cl''}$$
in $M$ then $x_{\cl''}=0$ for all $\cl''\preceq\cl$. We argue by induction on $\dim[\cl]-\dim[\cl'']$. If 
$\dim[\cl]-\dim[\cl'']=0$ then $\cl''=\cl$ and the equality $x_{\cl''}=0$ holds by assumption. Now assume that
$\dim[\cl]-\dim[\cl'']>0$ and that the result is known when $\cl''$ is replaced by $\cl'\in\fD^\s$ such that 
$\dim[\cl]-\dim[\cl']<\dim[\cl]-\dim[\cl'']$. We have
$$\sum_{\cl'\in\fD^\s;\cl''\preceq\cl'\preceq\cl}x_{\cl'}(u\i)\r_{\cl'',\cl'}=u^{-\dim[\cl]}x_{\cl''}(u).$$
Using the induction hypothesis, this becomes
$$x_{\cl''}(u\i)u^{\dim[\cl]-\dim[\cl'']}=x_{\cl''}(u).$$
This equality together with the condition that $\deg x_{\cl''}\le(\dim[\cl]-\dim[\cl'']-1)/2$ implies that
$x_{\cl''}=0$. This completes the inductive proof.

\subhead 5.3\endsubhead
We now specialize Theorem 5.2 in the context considered in 0.2. We obtain the following result.

\proclaim{Theorem 5.4} Let $w\in W^\s$. The polynomials $P^\s_{y,w}$ ($y\in W^\s$) in 0.2 are characterized by 
the following properties:
$$\sum_{y\in W^\s}P^\s_{y,w}(u\i)\ov{T_y}=u^{-\ul(w)}\sum_{y\in W^\s}P^\s_{y,w}(u)T_y\text{ in }\HH;\tag a$$
$$P^\s_{y,w}=0\text{ if }y\not\le w;\tag b$$
$$\deg P^\s_{y,w}\le(\ul(w)-\ul(y)-1)/2\text{ if }y'<w\text{ and }P^\s_{w,w}=1.\tag c$$
\endproclaim
This provides a new proof of the statement in \cite{\LC, (8.1)}. (Another proof, generalizing that in \cite{\KL},
was given earlier in \cite{\DG}.)

\head 6. A model for the $\HH$-module $M$\endhead
\subhead 6.1\endsubhead
We want to give a model of the $\HH$-module $M$ using functions on the set of rational points of a variety over 
a finite field analogous to the model of $\HH$ given in 4.4, 4.5. 

We preserve the setup of 2.1. Let $s$ be an odd integer $\ge1$. Note that $\tph^s\colon \cb@>>>\cb$ is the Frobenius 
map for an $\FF_{q^s}$-rational structure on $\cb$. Moreover $K$ is stable under $\tph^s\colon G@>>>G$. Let $\cf_s$ 
be the $\bbq$-vector space consisting of all functions $f\colon \cb^{\tph^s}@>>>\bbq$ which are constant on each 
$K^{\tph^s}$-orbit on $\cb^{\tph^s}$. For any $(\cs,\Ph)\in\cc$ we define an isomorphism 
$\Ph_s\colon \tph^{s*}\cs@>>>\cs$ as the composition
$$\tph^{s*}\cs@>\tph^{(s-1)*}(\Ph)>>\tph^{(s-1)*}\cs@>>>\do@>>>\tph^*\cs@>\Ph>>\cs$$
For any $B\in\cb^{\tph^s}$, $\Ph_s$ induces a linear isomorphism $\cs_B@>>>\cs_B$ whose trace is denoted by 
$\c_{s;\cs,\Ph}(B)$. For any $k\in K^{\tph^s}$ we have $\c_{s;\cs,\Ph}(kBk\i)=\c_{s;\cs,\Ph}(B)$. Hence the 
function $\c_{s;\cs,\Ph}$, $B\m\c_{s;\cs,\Ph}(B)$ belongs to $\cf_s$. Note that $(\cs,\Ph)\m\c_{s;\cs,\Ph}$ 
defines a group homomorphism $h\colon \frak K(\cc)@>>>\cf_s$. From the definitions we see that the kernel of $h$ contains 
$\Th(\frak K(\cc'))$ hence $h$ induces a group homomorphism $\vt_s\colon M@>>>\cf_s$. Note that $\vt_s$ induces a
$\bbq$-linear map
$$\bar\vt_s\colon \bbq\ot_\ca M@>>>\cf_s\tag a$$ 
where $\bbq$ is regarded as an $\ca$-algebra by $u\m q^s$. We have the following result.

(b) {\it If $K$ is connected then $\bar\vt_s$ is an isomorphism.}
\nl
The proof is the same as that of (24.2.7) in \cite{\CS} which deals with the $G$-action on the unipotent variety 
instead of the $K$-action on $\cb$.

Thus, in the setup of (a), $\bbq\ot_\ca M$ has, in addition to the basis 
$$\{(\cl,\b^\cl);\cl\in\fD^\s\}$$ 
(which depends on the choices of the $\a^\cl$ in 0.1), another basis
corresponding under $\bar\vt_s$ to the basis of $\cf_s$ given by the
characteristic functions of the $K^{\tph^s}$-orbits in  
$\cb^{\tph^s}$ (which is independent of any choice). Each element
$(\cl,\b^\cl)$ in the first basis is a linear  
combination with coefficients $\pm1$ of elements in the second basis
given by the characteristic functions of  
those $K^{\tph^s}$-orbits in $\cb^{\tph^s}$ which are contained in $[\cl]$.

\subhead 6.2\endsubhead
For $f\in\cg_s$ (see 4.4) and $f'\in\cf_s$ we define $f*f'\in\cf_s$ by
$$(f*f')(B)=\sum_{B'\in\cb^{\tph^s}}f(B,B')f'(B').$$
This defines a $\cg_s$-module structure on the $\bbq$-vector space $\cf_s$. Now let $(\fS,\Ph)\in\fC$, 
$(\cs',\Ph')\in\cc$. Using the definitions and Grothendieck's
``faisceaux-fonctions" dictionary we see that 
$$\vt_s(\fS,\Ph)*\vt_s(\cs',\Ph')=\vt_s((\fS,\Ph)*(\cs',\Ph'))$$ where
$*$ in the left hand side is as above, $*$ in the right hand side is
as in 4.1; the $\vt_s$ in the last equality are as in 4.5, 6.1, 6.1
respectively. It follows that, under the identification
$\bbq\ot_\ca\HH=\cg_s$ (see 4.5(a)) 

(a) {\it the linear map 6.1(a) is $\bbq\ot_\ca\HH$-linear.}
\nl
Note also that the following analogue of Lemma 4.6 holds: 

(b) {\it If $K$ is connected then the map $\vt\colon M@>>>\op_{s\in\{1,3,5,\do\}}\cf_s$ (with components $\vt_s$) is 
injective.}
\nl
This follows easily from 6.1(b). 

\subhead 6.3\endsubhead 
We would like to find an analogue of 6.2(b) without assuming that $K$
is connected. 

Let $x\in K$ and let $s$ be as in 6.1. We define a map  
$$\tph_{x,s}\colon G@>>>G, \qquad \tph_{x,s}(g)=x\tph^s(g)x\i.$$
Note that $\tph_{x,s}$ is the Frobenius map for an  
$\FF_{q^s}$-rational structure on $G$ (indeed, we have
$\tph_{x,s}=\Ad(y)\i\tph^s\Ad(y)$ where $y\in G$ is such  
that $x=y\i\tph^s(y)$). Moreover $K$ is stable under $\tph_{x,s}$. Now
$\tph_{x,s}$ induces a map $\cb@>>>\cb$  
denoted again by $\tph_{x,s}$ (it is the Frobenius map for an
$\FF_{q^s}$-rational structure on $\cb$). 
Let $\cf_{x,s}$ be the $\bbq$-vector space 
consisting of all functions $f\colon \cb^{\tph_{x,s}}@>>>\bbq$ which are constant on each $K^{\tph_{x,s}}$-orbit on 
$\cb^{\tph_{x,s}}$. For any $(\cs,\Ph)\in\cc$ we define an isomorphism 
$\Ph_{x,s}\colon \tph_{x,s}^*\cs@>>>\cs$ as the composition
$$\tph_{x,s}^*\cs=\tph^{s*}\Ad(x)^*\cs@>\tph^{s*}\ct_{x\i}>>\tph^{s*}\cs@>\tph^{(s-1)*}(\Ph)>>\tph^{(s-1)*}\cs@>>>\do@>>>\tph^*\cs@>\Ph>>\cs$$
where $\ct_{x\i}\colon \Ad(x)^*\cs@>>>\cs$ is given by the $K$-equivariant structure of $\cs$. For any 
$B\in\cb^{\tph_{x,s}}$, $\Ph_{x,s}$ induces a linear isomorphism $\cs_B@>>>\cs_B$ whose trace is denoted by 
$\c^{x,s}_{\cs,\Ph}(B)$. For any $k\in K^{\tph_{x,s}}$, $B\in\cb^{\tph_{x,s}}$ 
 we have $\c^{x,s}_{\cs,\Ph}(kBk\i)=\c^{x,s}_{\cs,\Ph}(B)$. 
Hence the function $\c^{x,s}_{\cs,\Ph}$, $B\m\c^{x,s}_{\cs,\Ph}(B)$ belongs to $\cf_{x,s}$. Note that 
$(\cs,\Ph)\m\c^{x,s}_{\cs,\Ph}$ defines a group homomorphism $h\colon \fK(\cc)@>>>\cf_{x,s}$. From the definitions we see 
that the kernel of $h$ contains $\Th(\fK(\cc'))$ hence $h$ induces a group homomorphism $\vt_{x,s}\colon M@>>>\cf_{x,s}$.
This induces a $\bbq$-linear map
$$\bar\vt_{x,s}\colon \bbq\ot_\ca M@>>>\cf_{x,s}\tag a$$
where $\bbq$ is regarded as an $\ca$-algebra by $u\m q^s$. 
It is likely that 

(b) {\it the map $\bbq\ot_\ca M@>>>\op_{x\in K}\cf_{x,s}$ (with components $\bar
\vt_{x,s}$) is injective.}
\nl
When $K$ is connected this follows from 6.1(b). (Note that $\vt_{1,s}=\vt_s$.)

\subhead 6.4\endsubhead
In the setup of 6.3 let $\utph_{x,s}=\tph_{x,s}\T\tph_{x,s}\colon \cb\T\cb@>>>\cb\T\cb$.
This is the Frobenius map for an $\FF_{q^s}$-rational structure on $\cb\T\cb$. Now $G^{\tph_{x,s}}$ acts on
$(\cb\T\cb)^{\utph_{x,s}}$ (by restriction of the $G$-action on $\cb\T\cb$) and from Lang's theorem we see that 
the $G^{\tph_{x,s}}$-orbits on $(\cb\T\cb)^{\utph_{x,s}}$ are exactly the sets $\fO_w^{\utph_{x,s}}$ where $w$ 
runs through $W^\s$.

Let $\cg_{x,s}$ be the $\bbq$-vector space consisting of all functions $f\colon (\cb\T\cb)^{\utph_{x,s}}@>>>\bbq$ which
are constant on each $G^{\tph_{x,s}}$-orbit on $(\cb\T\cb)^{\utph_{x,s}}$. For any $w\in W^\s$ let $f_{w,x,s}$ be
the function on $(\cb\T\cb)^{\utph_{x,s}}$ which is equal to $1$ on $\fO_w^{\utph_{x,s}}$ and is equal to $0$ on 
$\fO_{w'}^{\utph_{x,s}}$ for $w'\in W^\s-\{w\}$. Note that $\{f_{w,x,s};w\in W^\s\}$ is a $\bbq$-basis of 
$\cg_{x,s}$.

For $f,f'\in\cg_{x,s}$ we define $f*f'\in\cg_{x,s}$ by
$$(f*f')(B,B'')=\sum_{B'\in\cb^{\tph_{x,s}}}f(B,B')f'(B',B'').$$
This defines on $\cg_{x,s}$ a structure of associative $\bbq$-algebra with unit element $f_{1,x,s}$.
For $f\in\cg_{x,s}$ and $f'\in\cf_{x,s}$ we define $f*f'\in\cf_{x,s}$ by
$$(f*f')(B)=\sum_{B'\in\cb^{\tph_{x,s}}}f(B,B')f'(B').$$
This defines a $\cg_{x,s}$-module structure on the $\bbq$-vector space $\cf_{x,s}$. 

From \cite{\IW}, \cite{\MA} it follows that the following relations hold in the algebra $\cg_{x,s}$.

(a) $f_{w,x,s}*f_{w',x,s}=f_{ww',x,s}$ if $w,w'\in W^\s$, $\ul(ww')=\ul(w)+\ul(w')$;

(b) $f_{w_\o,x,s}*f_{w_\o,x,s}=q^{s\ul(w_\o)}f_{1,x,s}+(q^{s\ul(w_\o)}-1)f_{w_\o,x,s}$ if $\o\in\bS.$ 
\nl
Hence as in 4.5(a) we can identify $\bbq\ot_\ca\HH=\cg_{x,s}$ as $\bbq$-algebras by
$T_w\m f_{w,x,s}$ where $\bbq$ is regarded as an $\ca$-algebra with $u\m q^s$.
As in 6.2 we see that (under the identification $\bbq\ot_\ca\HH=\cg_{x,s}$),

(a) {\it the linear map 6.3(b) is $\bbq\ot_\ca\HH$-linear.}
\nl
It is likely that 

(b) {\it the map $\ti{\vt}\colon M@>>>\op_{x\in K,s\in\{1,3,5,\do\}}\cf_{x,s}$
(with components $\vt_{x,s}$) is injective.}
\nl
This would imply that the $\HH$-module structure on $M$ can be completely recovered from the $\HH$-module structures on the various $\cf_{x,s}$. 

When $K$ is connected, $\ti{\vt}$ is indeed injective, by 6.2(b). In the general case, (b) would follow from 6.3(b).

\head 7. Formulas for the action of ${\bold H}$ on $M$\endhead
\subhead 7.1. {\it Notation}\endsubhead
Recall from \S 0.2 that we write
$$(W,S), \qquad \sigma\colon S \rightarrow S\eqno(7.1)(a)$$
for the Weyl group of $G$ and its automorphism induced by the
automorphism of $G$.  We also write
$$\overline{S} = \hbox{orbits of $\sigma$ on $S$}.\eqno(7.1)(b)$$
For each orbit $\omega \in {\overline S}$, we write
$$w_\omega = \text{long element of subgroup $W(\omega)$ generated by
$\omega$}. \eqno(7.1)(c)$$
Because $\sigma$ is an involutive automorphism, there are three
possibilities for $W(\omega)$:
$$W(\omega) = \cases S_2, & \omega= \{s\} \subset S; \\
S_2 \times S_2, & \omega = \{s,t\} \subset S,\  st=ts;\\
S_3, & \omega = \{s,t\} \subset S, \ (st)^3 = 1. \endcases \eqno(7.1(d)$$
We call these three cases {\it types one, two, and three}; therefore
$$\underline{\ell}(w_\omega) = m \text{\ if $\omega$ is type
$m$}. \eqno(7.1)(e)$$ 
Here $\underline{\ell}$ is the length function on $W$ (\S 0.2).
We will abuse notation and identify
$$\overline{S} \simeq \{w_\omega\} \subset W^\sigma, \eqno (7.1)(f),$$
allowing us to write
$$(W^\sigma,\overline{S}) \eqno(7.1)(g)$$ 
for the Coxeter group presentation of $W^\sigma$. (The group
$W^\sigma$ is in fact the Weyl group of the reductive subgroup of $G$
fixed by a distinguished (that is, preserving some pinning)
automorphism $\sigma_d$ inner to $\sigma$, but we will make no use of
this fact.)

\subhead 7.2\endsubhead
In this section we will make explicit the action of the generators
$$T_{w_\omega},\qquad (T_{w_\omega} + 1)(T_{w_\omega} - u^m) = 0
\qquad (\omega \in \overline{S} \ \text{type $m$}) \eqno (7.2)(a)$$ 
of the quasisplit (unequal parameter) Iwahori Hecke
algebra ${\bold H}$ (introduced geometrically in \S 3.1, with the
algebra structure defined in \S 4.2, and identified with the Iwahori
Hecke algebra in \S 4.7)  on the basis
$$a_{\Cal L} = \{({\Cal L},\beta^{\Cal L}) \mid {\Cal L} \in {\frak
D}^\sigma \} \eqno(7.2)(b)$$ 
for the module $M$ introduced in \S 2.3 (as an ${\Cal A}$-module)
and \S 4.2 (as an ${\bold H}$-module). In order to simplify the
notation slightly we will in this section write simply ${\Cal L}$
instead of $a_{\Cal L}$; this should cause no ambiguity or confusion.

What we are going to see is that the matrix of $T_{w_\omega}$ is
block-diagonal, with blocks of size one, two, three, or four. (The
corresponding partition of the basis elements ${\Cal L}$ is different
for each $w_\omega$, so the whole action of ${\bold H}$ need {\it not}
be block-diagonal.)

At least over the quotient field of ${\Cal A}$, the quadratic relation
$(7.2)(a)$ guarantees that $M$ is the direct sum of
$$\aligned \text{$u^m$-eigenspace of $T_{w_\omega}$} &= \text{image of
$T_{w_\omega}+1$}, \\ \text{$-1$-eigenspace of $T_{w_\omega}$} &=
\text{kernel of $T_{w_\omega} +1$} \endaligned \eqno(7.2)(c)$$
It will be useful (for the recursion algorithm described in \S 8) to
write these eigenspaces also as we go along. 

The operator $T_{w_\omega} + 1$ is described geometrically in
(4.8). Write ${\Cal P}_{w_\omega}$ for the partial flag variety of
parabolic subgroups of type $\omega$, and 
$$\pi_{\omega}\colon {\Cal B}\rightarrow {\Cal P}_{w_\omega}$$
for the corresponding projection. The (same) two projections of ${\Cal
B}$ onto ${\Cal P}_{w_\omega}$ form a Cartesian square with the
projections $\pi_i\colon \overline{{\frak O}_{w_\omega}} \rightarrow
{\Cal B}$ considered in (4.8)(c). Proper base change and (4.8)(d)
imply that $T_{w_\omega} + 1$ is implemented on the level of sheaves
by $\pi_\omega^*\pi_{\omega!}$. Roughly speaking, it follows that
$$\aligned \text{image of ${T_{w_\omega} + 1}$ } &\longleftrightarrow
\text{sheaves 
pulled back from ${\Cal P}_{w_\omega}$;}\\ \text{kernel of
${T_{w_\omega} + 1}$ } &\longleftrightarrow \text{sheaves pushing forward
to zero on ${\Cal P}_{w_\omega}$.} \endaligned\eqno(7.2)(d)$$

For the intersection homology complexes ${\Cal L}^\sharp$ in which we
are ultimately interested, and the corresponding basis elements
${\frak A}_{\Cal L}$ of $M$, it turns out that ${\frak A}_{\Cal L}$
belongs to the $u^m$-eigenspace of $T_{w_\omega}+1$ if and only if
${\Cal L}^\sh$ is pulled back from an intersection cohomology complex
on ${\Cal P}_{w_\omega}$. We call these $w_\omega$ the {\it descents
for ${\Cal L}$}, by analogy with the corresponding behavior for
Schubert cells and $W$.  The remaining cases are called {\it ascents
for ${\Cal L}$}. We use this terminology to help sort the cases below,
without explicitly verifying the corresponding geometric properties.

\subhead 7.3\endsubhead
It will be useful to consider also the split Iwahori Hecke algebra
$${\Cal H} = {\Cal H}(W) = \langle T_s \mid s\in S\rangle, \eqno (7.3)(a)$$
with basis $\{T_w\mid w\in W\}$ as an $\Cal A$-module.  This Hecke
algebra has a module
$${\Cal M} = \sum_{{\Cal L} \in {\frak D}} {\Cal A} \cdot {\Cal L}
\eqno (7.3)(b)$$
defined in \cite{LV1, \S 3.1}. We have already noted in \S 4.9 a
relationship between the action of ${\bold H}$ on $M$ and that of
${\Cal H}$ on ${\Cal M}$, using the obvious forgetful map
$${\frak D}^\sigma \hookrightarrow {\frak D}$$
on parameters.  In order to recall that action of ${\Cal H}$, we need
to recall an explicit parametrization of ${\frak D}$ going back to
Kostant, Wolf, and Matsuki. 
\proclaim{Theorem 7.4} 

(i) There are in $G$ finitely many orbits $\{{\Cal S}_1,\ldots,{\Cal
S}_r\}$ under $K$ of $\theta$-stable maximal tori. Each orbit has a
representative $H_i$ preserved by the split Frobenius morphism $\phi$
of \S 1.1. 

(ii) The orbits of (i) are
permuted by $\sigma$.  The representative $H_i$ of each $\sigma$-fixed
orbit ${\Cal S}_i$ may be chosen to be preserved by $\sigma$, and
therefore by the Frobenius automorphism $\widetilde \phi$ of \S 2.1.

(iii) Every Borel subgroup $B$ of $G$ contains a $\theta$-stable maximal
torus $H$, unique up to conjugation by $B\cap K$. Consequently every
orbit of $K$ on $\Cal B$ has a representative $B_j$ containing one of
the representative $\theta$-stable tori $H_{i(j)}$.

(iv) In the setting of (ii), the $K$-equivariant local systems on the
orbit $K\cdot B_j$ are naturally parametrized by the characters of the
component group $H_{i(j)}^\theta/(H_{i(j)}^\theta)_0$, which is an
elementary abelian $2$-group. This character group can be naturally
described in terms of the lattice $X^*(H_{i(j)})$ of rational characters:
$$[H_{i(j)}^\theta/(H_{i(j)}^\theta)_0]\,\,\widehat{} \ = \
X^*(H_{i(j)})^{-\theta}/(1-\theta)X^*(H_{i(j)}).$$

(v) Write $W(G,H_i) = N_G(H_i)/H_i$ for the Weyl group. Then
$W(G,H_i)$ acts in a simply transitive fashion on the set of Borel
subgroups 
$${\Cal B}^{H_i} = \{B\in {\Cal B} \mid H_i \subset B\};$$
fixing one of these defines an isomorphism
$$ i_B\colon W(G,H_i)\rightarrow W.$$ 
Write
$$W(K,H_i) = N_K(H_i)/H_i \cap K$$
for the subgroup having representatives in $K$. Then the 
orbits of $K$ on ${\Cal B}$ corresponding to $H_i$ are in one-to-one
correspondence with the orbits of $W(K,H_i)$ on ${\Cal B}^{H_i}$. The
number of such orbits is therefore equal to the index of $W(K,H_i)$ in
$W(G,H_i)$. 
\endproclaim
\subhead 7.5\endsubhead Theorem 7.4 provides some additional
structure on the set $S$ of
generators of $W$ attached to a parameter 
$${\Cal L} \in {\frak D}.\eqno(7.5)(a)$$
(In certain cases the parameter will be naturally one element of a pair; in
those cases we will write ${\Cal L}_1$ instead of ${\Cal L}$, and
${\Cal L}_2$ for the other element of the pair.) To see this
structure, first write 
$$\ell({\Cal L}) = \dim([{\Cal L}]), \eqno(7.5)(b)$$
the dimension of the underlying $K$ orbit on ${\Cal B}$.
Now fix a representative $B\supset H_i$ of
$[{\Cal L}]$ as in the theorem. The roots of $H_i$ in $B$ define a
system of positive roots
$$R^+ =_{\text{def}} R(B,H_i) \subset R(G,H_i). \eqno(7.5)(b')$$
Write
$$\chi = \chi_{{\Cal L},B,H_i} \in
[H_i^\theta/(H_i^\theta)_0]\,\,\widehat{} \eqno(7.5)(b'')$$
for the character of the component group corresponding to ${\Cal L}$.

Use the isomorphism $i_B$ of Theorem 7.4(v) to
identify $W$ with $W(G,H_i)$, and therefore $s\in S$ with a simple
root $\alpha\in R^+$. For $w\in W$, we define (following
\cite{V, Definition 8.3.1})
$$w\times {\Cal L} = (K\cdot (w^{-1}\cdot B),{\Cal L}_w); \eqno(7.5)(c)$$
here we act on $B\in {\Cal B}^{H_i}$ using the isomorphism $i_B$ of
$W$ with $W(G,H_i)$, and we use the line bundle ${\Cal L}_w$
corresponding to the character 
$$\chi_w = \chi + \sum\Sb \alpha \in R^+\\
\theta\alpha=-\alpha\\ w\alpha 
\notin R^+\endSb \alpha. \eqno(7.5)(c')$$ 
of $H_i^\theta/(H_i^\theta)_0$. (Theorem 7.4(iv) shows how to
interpret the roots in the sum---called {\it real} below---as
characters of the component group.)

Because $H_i$ is preserved by the automorphism $\theta$ of $G$,
$\theta$ induces an automorphism of the root system $R(G,H_i)$.
We say that $s$ is {\it complex for ${\Cal L}$} if $\theta\alpha \ne \pm
\alpha$; equivalently, if $\ell(s\times {\Cal L}) = \ell({\Cal L}) \pm
1$. This is in some sense the most common situation. 

We say that
$s$ is a {\it complex ascent for ${\Cal L}$} 
$${\Cal L}' =_{\text{def}} s\times{\Cal L} \text{\ \ satisfies\ \ }
\ell(s\times {\Cal L}) = \ell({\Cal L}) + 1. \eqno(7.5)(d)$$ 
It is equivalent to require
$$\theta\alpha \in R^+\backslash\{\alpha\}. \eqno(7.5)(d')$$ 
In this case the action of the generator for ${\Cal H}$ is (\cite{LV1,
Lemma 3.5(b)})
$$T_s\cdot {\Cal L} = {\Cal L}', \qquad T_s\cdot{\Cal L}' = u{\Cal L}
+ (u-1){\Cal L}' \eqno(7.5)(d'')$$

Similarly, we say that $s$ is a {\it complex descent for ${\Cal L}'$};
this is characterized by 
$${\Cal L} =_{\text{def}} s\times{\Cal L}' \text{\ \ satisfies\ \ }
\ell({\Cal L}) = \ell({\Cal L}') - 1. \eqno(7.5)(e)$$ 

The eigenspaces of $T_s$ on the span of ${\Cal L}$ and ${\Cal L}'$ are
$$\text{$u$-eigenspace} = \langle {\Cal L} + {\Cal L}'\rangle, \quad
\text{$-1$-eigenspace} = \langle u{\Cal L} - {\Cal L}'\rangle. \eqno(7.5)(e')$$

We say that $s$ (or the root $\alpha$) is {\it imaginary for ${\Cal L}$}
if $\theta\alpha = 
\alpha$. In this case $\theta$ must preserve the root space ${\frak
g}_\alpha$; we say that 
$$\gathered \text{$s$ (or $\alpha$) is {\it noncompact imaginary for
${\Cal L}$} 
if $\theta|_{{\frak g}_{\alpha}} = -1$}\\ \text{$s$ (or $\alpha$) is
{\it compact imaginary for ${\Cal L}$} 
if $\theta|_{{\frak g}_{\alpha}} = +1$}\endgathered \eqno(7.5)(f)$$
Similarly, we say that
$s$ (or $\alpha$) is {\it real for ${\Cal L}$} if $\theta\alpha =
-\alpha$. In this case
$$m_\alpha =_{\text{def}} \alpha^\vee(-1) \in H^\theta = H\cap K.$$
We say that 
$$\gathered \text{$s$ (or the coroot $\alpha^\vee$) {\it satisfies the parity
condition for ${\Cal L}$} if $\chi(m_\alpha) = 1$;}\\ \text{$s$ {\it
does not satisfy the parity condition for ${\Cal L}$} if
$\chi(m_\alpha) = -1$.}\endgathered \eqno(7.5)(g)$$ 
Here $\chi$ is the character of $(H\cap K)/(H\cap K)_0$ from
$(7.5(b'')$.

We say that $s$ is an {\it imaginary noncompact type I ascent for
${\Cal L}_1$} if it is noncompact, and $s_\alpha \notin W(K,H_i)$; the
second condition is equivalent to 
$${\Cal L}_2 =_{\text{def}} s \times {\Cal L}_1 \ne {\Cal
L}_1. \eqno(7.5)(h)$$  
In this case the construction of \cite{LV1, Lemma 3.5(d1)} defines a
{\it single-valued Cayley transform}
$$c_s({\Cal L}_1) = c_s({\Cal L}_2) = \{{\Cal L}'\}, \qquad \ell({\Cal L}') =
\ell({\Cal L}_j)+1. \eqno(7.5)(h')$$
The corresponding $\theta$-stable maximal torus $c_s(H_i)$ can be
constructed using the root $SL(2)$ for $\alpha$.  The formulas for the
Hecke algebra action are 
$$\gathered T_s\cdot {\Cal L}_1 = {\Cal L}_2 + {\Cal L}', \qquad T_s\cdot {\Cal
L}_2 = {\Cal L}_1 + {\Cal L}',\\
T_s\cdot {\Cal L}' = (u-2){\Cal L}' + (u-1)({\Cal L}_1 + {\Cal
L}_2). \endgathered\eqno(7.5)(h'')$$ 

We say that $s$ is a {\it real type I descent} for ${\Cal
L}'$;  it is equivalent to say that $s$ is real, satisfies the parity
condition for ${\Cal L}'$, and that the root 
$\alpha$ does not takes the value $-1$ on $c_s(H_i)\cap K$; equivalently, if
$s \times {\Cal L}' = {\Cal L}'$. In this case ${\Cal L}'$ is equal to
the Cayley transform through $s$ of exactly ${\Cal L}_1$ and ${\Cal
L}_2$: we write 
$$c^s({\Cal L}') =_{\text{def}} \{{\Cal L}_1,{\Cal L}_2\},
\eqno(7.5)(i)$$
and call this two-element set a {\it double-valued inverse Cayley
transform}. 

The eigenspaces of $T_s$ on the span of ${\Cal L}_1$, ${\Cal L}_2$, and
${\Cal L}'$ are 
$$\gathered\text{$u$-eigenspace} = \langle {\Cal L}_1 + {\Cal L}_2 + {\Cal
L}'\rangle,\\ \text{$-1$-eigenspace} =  \langle(u-1){\Cal L}_1 -
{\Cal L}', (u-1){\Cal L}_2 -{\Cal L}'\rangle.\endgathered \eqno(7.5)(i')$$

We say that $s$ is an {\it imaginary noncompact type II ascent for
${\Cal L}$} if it is noncompact, and $s_\alpha \in W(K,H_i)$; the
second condition is equivalent to 
$$s \times {\Cal L} = {\Cal L}. \eqno(7.5)(j)$$ 
In this case the construction of \cite{LV1, Lemma 3.5(c1)} defines a
{\it double-valued Cayley transform}
$$c_s({\Cal L}) = \{{\Cal L}'_1,{\Cal L}'_2\},\quad {\Cal L}'_2 =
s\times {\Cal L}'_1, \quad \ell({\Cal L}'_i) =
\ell({\Cal L})+1. \eqno(7.5)(j')$$
The formulas for the Hecke algebra action are
$$\gathered T_s\cdot {\Cal L} = {\Cal L} + ({\Cal L}'_1 + {\Cal L}'_2),
\\ T_s\cdot {\Cal L}'_1 = (u-1){\Cal L} + (u-1){\Cal L}'_1 - {\Cal
L}'_2,\\  T_s\cdot {\Cal L}'_2 = (u-1){\Cal L} - {\Cal L}'_1
+(u-1){\Cal L}'_2.\endgathered \eqno(7.5)(j'')$$ 

We say that $s$ is a {\it real type II descent for ${\Cal
L}'_j$}; it is equivalent to say that $\alpha$ is real, satisfies the parity condition for ${\Cal
L}'_j$, and takes the value $-1$ on $c_s(H_i)\cap K$.
In this case the pair $\{{\Cal L}'_1,{\Cal L}'_2\}$ is the Cayley
transform through $s$ of exactly one local system:
we write
$$c^s({\Cal L}'_1) = c^s({\Cal L}'_2) = \{{\Cal L}\}, \eqno(7.5)(k)$$
which we call {\it single-valued inverse Cayley
transforms}. 

The eigenspaces of $T_s$ on the span of ${\Cal L}$, ${\Cal L}'_1$, and
${\Cal L}'_2$ are 
$$\gathered \text{$u$-eigenspace} = \langle {\Cal L}+ {\Cal L}'_1, {\Cal L} +
{\Cal L}'_2\rangle, \\ \text{$-1$-eigenspace} =  \langle(u-1){\Cal L} -
{\Cal L}'_1 -{\Cal L}'_2\rangle. \endgathered \eqno(7.5)(k')$$

We say that $s$ is a {\it real nonparity ascent for ${\Cal L}$} if
$\alpha$ is real and fails to satisfy the parity condition for ${\Cal
L}$ (cf. $(7.5)(g')$). In this case (\cite{LV1, Lemma 3.4(e)})
$$T_s({\Cal L}) =  -{\Cal L}, \quad \text{$-1$-eigenspace} =
\langle{\Cal L}\rangle, \quad \text{$u$-eigenspace} = 0.\eqno(7.5)(\ell)$$

Finally, we say that $s$ is a {\it compact imaginary descent for
${\Cal L}$} if $\alpha$ is compact imaginary. In this case
(\cite{LV1, Lemma 3.4(a)})
$$T_s({\Cal L}) =  u{\Cal L}, \quad \text{$-1$-eigenspace} =  0, \quad
\text{$u$-eigenspace} = \langle {\Cal L}\rangle.\eqno(7.5)(m)$$

If $s \in S$ is a $\sigma$-fixed generator, and ${\Cal L} \in {\frak
D}^\sigma$, then the formulas (7.5) essentially give the action of
$T_s \in {\bold H}$. There are two exceptions. First, if $s$ is a real
type I descent for ${\Cal L}'$, it may happen that $\sigma$
interchanges the two elements $\{{\Cal L}_1,{\Cal L}_2\}$ of
$c^s({\Cal L}')$. In this case the formulas $(7.5)(h'')$ must be
replaced by the single formula
$$T_s({\Cal L}') = u{\Cal L}' \qquad ({\Cal L}_i \notin {\frak
D}^\sigma). \eqno(7.5)(h''')$$
Similarly, if $s$ is a type II ascent for ${\Cal L}$, it may happen
that $\sigma$ interchanges the two elements $\{{\Cal L}'_1,{\Cal L}'_2\}$ of
$c_s({\Cal L})$. In this case the formulas $(7.5)(j'')$ are
replaced by 
$$T_s({\Cal L}) = -{\Cal L} \qquad ({\Cal L}'_i \notin {\frak
D}^\sigma). \eqno(7.5)(j''')$$

\subhead 7.6. Formulas in the type 2 cases \endsubhead Now suppose 
$${\Cal L} \in {\frak D}^\sigma,$$
and $(H_i,B)$ is chosen as in (7.5). Write 
$$R^+ = R(B,H_i) \subset R(G,H_i)$$
for the corresponding positive root system.
Fix a pair $\omega = \{s,t\}$ of commuting simple reflections interchanged by
the automorphism $\sigma$, so that 
$$w_\omega = s t \in W^\sigma$$
is a simple generator of length two.  The assumption that ${\Cal L}$ is
fixed by $\sigma$ means that 
$$(\sigma(B),\sigma(H_i)) = \text{Ad}(k)(B,H_i)$$
for some $k\in K$; the coset $k(B\cap K)$ is uniquely determined by
Theorem 7.4(iii). Consequently $\text{Ad}(k^{-1})\sigma$ defines an
automorphism, which we will simply write $\sigma$, of the based root
datum corresponding to $(B,H_i)$. This description makes it clear that
the automorphism $\sigma$ of $W(G,H_i) \simeq W$ (Theorem 7.4) must
preserve the subgroup $W(K,H_i)$, as well as the ``status'' (complex
ascent, noncompact imaginary, etc.) of the simple roots described in
(7.5). In particular, the simple roots $\alpha$ and $\beta$ of $R^+$
corresponding to the reflections $s$ and $t$ must have the same status.

We will enumerate twelve cases, and in each case give a formula for
the action of the generator $T_{w_\omega}$ on ${\Cal L}$. 

We say that $w_\omega$ is a {\it two-complex ascent for ${\Cal L}$} if
$$\theta\alpha \in R^+, \qquad \theta\alpha \notin
\{\pm\alpha,\pm\beta\}. \eqno(7.6)(a)$$ 
Write
$${\Cal L}' = w_\omega \times {\Cal L} = st \times {\Cal L}, \qquad
\ell({\Cal L}')= \ell({\Cal L})+2;\eqno(7.6)(a')$$
the last length condition is equivalent to the definition of
two-complex ascent in(7.6)(a).  In this case
$$T_{w_\omega}({\Cal L}) = {\Cal L}', \quad  T_{w_\omega}({\Cal L}') =
u^2{\Cal L} + (u^2-1){\Cal L}'. \eqno(7.6)(a'')$$

We say that $w_\omega$ is a {\it two-complex
descent for ${\Cal L}'$}. 
These first two cases generalize \cite{LV2, Theorem 0.2(iii) and
(iv)}.  The eigenspaces of $T_{w_\omega}$ on the span of ${\Cal L}$
and ${\Cal L}'$ are 
$$\text{$u^2$-eigenspace} = \langle {\Cal L} + {\Cal L}'\rangle, \quad
\text{$-1$-eigenspace} = \langle u^2{\Cal L} - {\Cal
L}'\rangle. \eqno(7.6)(b)$$ 

We say that $w_\omega$ is a {\it two-semiimaginary ascent for
${\Cal L}$} if
$$\theta\alpha =\beta\in R^+. \eqno(7.6)(c)$$ 
Define the {\it Cayley transform of ${\Cal L}$ through $w_\omega$} to be
$$c_{w_\omega}({\Cal L}) =_{\text{def}}\{ s\times {\Cal L}\} = \{t \times
{\Cal L}\} =_{\text{def}} \{{\Cal L}'\}, \qquad
\ell({\Cal L}')= \ell({\Cal L})+1;\eqno(7.6)(c')$$
the second equality arises by applying $s_\alpha s_\beta \in W_K(H_i)$
to the first. In this case
$$T_{w_\omega}({\Cal L}) = u{\Cal L} + (u+1){\Cal L}', \quad
T_{w_\omega}({\Cal L}') =
(u^2-u){\Cal L} + (u^2-u-1){\Cal L}'. \eqno(7.6)(c'')$$
We say that $w_\omega$ is a {\it two-semireal
descent for ${\Cal L}'$}. The {\it inverse Cayley transform} is
$$c^{w_\omega}({\Cal L}') =_{\text{def}} \{s\times {\Cal L}'\} = \{t \times
{\Cal L}'\} = \{{\Cal L}\}. \eqno(7.6)(d)$$
These second two cases generalize \cite{LV2, Theorem 0.2(i) and
(ii)}.  The eigenspaces of $T_{w_\omega}$ on the span of ${\Cal L}$
and ${\Cal L}'$ are 
$$\text{$u^2$-eigenspace} = \langle {\Cal L} + {\Cal L}'\rangle, \quad
\text{$-1$-eigenspace} = \langle (u^2-u){\Cal L} - (u+1){\Cal
L}'\rangle. \eqno(7.5)(d')$$ 

These four cases exhaust the possibilities when $\alpha$
and $\beta$ are neither real nor imaginary.  Next we consider the real and
imaginary cases.  We say that $w_\omega$ is a {\it two-imaginary
noncompact type I-I ascent for ${\Cal L}_1$} if $\alpha$ and $\beta$
are noncompact 
imaginary roots in $R^+$, and $W(K,H_i) \cap \langle
s_\alpha,s_\beta \rangle$ has just one element.  This case arises in
the example $(9.1)(f)$ below, and the calculations we omit
can be carried out almost entirely in that example. We define
$${\Cal L}_2 = w_\omega\times {\Cal L}_1, \qquad {\Cal L}' =
c_t(c_s({\Cal L})),\qquad \ell({\Cal L}') = \ell({\Cal L}_j)+2;
\eqno(7.6)(e)$$  
the Cayley transforms are defined by the ``noncompact imaginary''
hypothesis, and single-valued by the ``type I-I'' hypothesis.  We have
$$\gathered T_{w_\omega}({\Cal L}_1) = {\Cal L}_2 + {\Cal L}', \quad
T_{w_\omega}({\Cal L}_2) = {\Cal L}_1 + {\Cal L}', \\
T_{w_\omega}({\Cal L}') = (u^2-1)({\Cal L}_1+{\Cal L}_2) +
(u^2-2){\Cal L}'. \endgathered\eqno(7.6)(e')$$

We say that $w_\omega$ is {two-real type I-I descent for ${\Cal L}'$}.
The eigenspaces of $T_{w_\omega}$ on the span of ${\Cal L}_1$, ${\Cal
L}_2$, and ${\Cal L}'$ are 
$$\gathered \text{$u^2$-eigenspace} = \langle {\Cal L}_1 + {\Cal L}_2 + {\Cal
L}'\rangle, \\
\text{$-1$-eigenspace} =  \langle(u^2-1){\Cal L}_1 -
{\Cal L}', (u^2-1){\Cal L}_2 -{\Cal L}'\rangle. \endgathered\eqno(7.6)(f')$$

We say that $w_\omega$ is a {\it two-imaginary noncompact type II-II
ascent} if $\alpha$ and $\beta$ are noncompact imaginary roots
in $R^+$, and $W(K,H_i) \cap \langle s_\alpha,s_\beta\rangle$ has
four elements.  (This case arises in the example $(9.1)(h)$ of the next
section.) Consequently (in fact equivalently) the imaginary 
roots $\alpha$ and $\beta$ are type II, so the two Cayley transforms
$c_s({\Cal L})$ and $c_t({\Cal L})$ 
are double-valued, and the iterated Cayley transform takes four
values:
$$\{{\Cal L}', s\times{\Cal L}', t\times {\Cal L}', st\times{\Cal L}'\}
= c_s(c_t({\Cal L}) = c_t(c_s({\Cal L})). \eqno(7.6)(g)$$
Two of these four elements, which we call ${\Cal L}'_1$ and ${\Cal
L}'_2$, belong to ${\frak D}^\sigma$, and the other 
two are interchanged by $\sigma$; we define the {\it double-valued
Cayley transform}
$$c_{w_\omega}({\Cal L}) = c_s(c_t({\Cal L})) \cap {\frak D}^\sigma =
\{{\Cal L}'_1,{\Cal L}'_2\}, \eqno(7.6)(g')$$
with the ${\Cal L}'_j$ interchanged by the cross action of $st=w_\omega$. Then
$$\gathered T_{w_\omega}({\Cal L}) = {\Cal L} + {\Cal L}'_1 +
{\Cal L}'_2, \\ 
T_{w_\omega}({\Cal L}'_1) = (u^2-1){\Cal L} + (u^2-1){\Cal L}'_1 -
{\Cal L}'_2,\\  
T_{w_\omega}({\Cal L}'_2) = (u^2-1){\Cal L} -{\Cal L}'_1 +(u^2-1)
{\Cal L}'_2 . \endgathered\eqno(7.6)(g'')$$  

We say that $w_\omega$ is a {\it two-real type II-II
descent for ${\Cal L}'_j$}, and define the inverse Cayley transform
$$c^{w_\omega}({\Cal L}'_j) = \{{\Cal L}\}.\eqno(7.6)(h)$$

The eigenspaces of $T_{w_\omega}$ on the span of ${\Cal L}$, ${\Cal
L}'_1$, and ${\Cal L}'_2$ are 
$$\gathered \text{$u^2$-eigenspace} = \langle {\Cal L}+ {\Cal L}'_1,
{\Cal L} + {\Cal L}'_2\rangle, \\
\text{$-1$-eigenspace} =  \langle(u^2-1){\Cal L} -
({\Cal L}'_1 +{\Cal L}'_2)\rangle. \endgathered\eqno(7.6)(h')$$

We say that $w_\omega$ is a {\it two-imaginary noncompact type I-II
ascent for ${\Cal L}_1$} if $\alpha$ and $\beta$ are noncompact imaginary roots
in $R^+$, and $W(K,H_i) \cap \langle s_\alpha,s_\beta\rangle$ has
two elements.  Because this last intersection is preserved by the
automorphism $\sigma$, it must consist of the identity and $s_\alpha
s_\beta$. (This case arises in the example $(9.1)(g)$ below.)  Define
$${\Cal L}_2 = s \times {\Cal L}_1 = t\times {\Cal L}_1 \in {\frak
D}^\sigma, \quad \ell({\Cal L}_2) = \ell({\Cal L}_1).\eqno(7.6)(i)$$
The imaginary roots $\alpha$ and $\beta$ are  
type I (because $s$ and $t$ do not belong to $W(K,H_i)$, so the two
Cayley transforms 
$$c_s({\Cal L}_1) = c_s(s\times{\Cal L}_1) = c_s({\Cal L}_2)$$ 
and $c_t({\Cal L}_1) = c_t({\Cal L}_2)$
are single-valued:
$$\sigma\cdot c_s({\Cal L}_j) = c_t({\Cal L}_j) \in {\frak D}.$$

Write $H_\beta$ for the $\theta$-stable Cartan underlying
$c_t({\Cal L}_1)$ (the Cayley transform of $H_i$ through the imaginary
root $\beta$.) The presence of $s_\alpha s_\beta$ in $W(K,H_i)$
implies that the 
simple reflection corresponding to $s$ belongs to
$W(K,H_\beta)$. Therefore the simple reflection $s$ is type II
noncompact imaginary for 
$c_t({\Cal L})$, so the iterated Cayley transform is double-valued. We define
$$c_{w_\omega}({\Cal L}_1) = c_{w_\omega}({\Cal L}_2) = c_s(c_t({\Cal
L})) = \{{\Cal L}'_1,{\Cal L}'_2\} \subset
{\frak D}^\sigma, \quad \ell({\Cal L}_j') = \ell({\Cal L}_1) +2. \eqno(7.6)(i')$$
Then we can choose the rational structures so that
$$\gathered T_{w_\omega}({\Cal L}_1) = {\Cal L}_1 + {\Cal L}'_1 + {\Cal L}'_2,
\qquad T_{w_\omega}({\Cal L}_2) = {\Cal L}_2 + {\Cal L}'_1 - {\Cal
L}'_2,\\   T_{w_\omega}({\Cal L}'_1) = (u^2-1)({\Cal L}_1+{\Cal L}_2) + 
(u^2-2){\Cal L}'_1, \\
T_{w_\omega}({\Cal L}'_2) = (u^2-1)({\Cal L}_1 - {\Cal L}_2)
+ (u^2-2){\Cal L}'_2\endgathered \eqno(7.6)(i'')$$  

We say that $w_\omega$ is a {\it two-real type II-I
descent for ${\Cal L}'_j$}, and define the double-valued inverse
Cayley transform by  
$$c^{w_\omega}({\Cal L}'_j) = \{{\Cal L}'_1,{\Cal L}'_2\} =
c^s(c^t({\Cal L}'_j)) \subset
{\frak D}^\sigma,\eqno(7.6)(j)$$

The eigenspaces of $T_{w_\omega}$ on the span of ${\Cal L}_1$, ${\Cal
L}_2$, ${\Cal L}'_1$, and ${\Cal L}'_2$ are 
$$\gathered \text{$u^2$-eigenspace} = \langle {\Cal L}_1+ {\Cal L}_2 + {\Cal
L}'_1, {\Cal L}_1 - {\Cal L}_2 + {\Cal L}'_2 \rangle, \\
\text{$-1$-eigenspace} =  \langle(u^2-1){\Cal L}_1 - ({\Cal L}'_1
+{\Cal L}'_2), (u^2-1){\Cal L}_2 - ({\Cal L}'_1 - {\Cal L}'_2)
\rangle. \endgathered \eqno(7.6)(j')$$ 

We say that $w_\omega$ is a {\it two-real nonparity ascent} if
$\alpha$ and $\beta$ are real roots not satisfying the parity
condition (\cite{LV1, Lemma 3.5(e)}).  In this case
$$T_{w_\omega}({\Cal L}) =  -{\Cal L}, \quad \text{$u^2$-eigenspace} =
0, \quad \text{$-1$-eigenspace} = \langle{\Cal L}\rangle.\eqno(7.6)(k)$$

Finally, we say that $w_\omega$ is a {\it two-imaginary compact
descent} if $\alpha$ and $\beta$ are compact imaginary roots. In
this case
$$T_{w_\omega}({\Cal L}) = u^2{\Cal L}, \quad \text{$u^2$-eigenspace}
= \langle{\Cal L}\rangle, \quad \text{$-1$-eigenspace} = 0.\eqno(7.6)(\ell)$$

\subhead 7.7. Formulas in the type 3 cases \endsubhead We retain the
notation introduced at the beginning of \S 7.6; but now suppose 
$\omega = \{s,t\}$ is a pair of {\it noncommuting} simple reflections
interchanged by the automorphism $\sigma$, so that 
$$w_\omega = s t s = tst \in W^\sigma$$
is a simple generator of length three.  Again we write $\alpha$ and
$\beta$ for the simple roots in $R^+$ (of $H_i$ in $B$) corresponding
to $s$ and $t$; again the automorphism $\sigma$ ensures that $\alpha$
and $\beta$ must have the same ``status'' as described in (7.5). We
will enumerate eight cases, and in each case give a formula for the
action of the generator $T_{w_\omega}$ on ${\Cal L}$. 

We say that $w_\omega$ is a {\it three-complex ascent for ${\Cal L}$} if
$$\theta\alpha \in R^+, \qquad \theta\alpha \notin
\{\pm\alpha,\,\pm\beta,\,\pm(\alpha + \beta)\}. \eqno(7.7)(a)$$ 
Write
$${\Cal L}' = w_\omega \times {\Cal L} = sts\times {\Cal L}, \qquad
\ell({\Cal L}')= \ell({\Cal L})+3;\eqno(7.7)(a')$$
the last length condition is equivalent to the definition of
three-complex ascent in (7.7)(a).  In this case
$$T_{w_\omega}({\Cal L}) = {\Cal L}', \qquad T_{w_\omega}({\Cal L}') =
u^3{\Cal L} + (u^3-1){\Cal L}'. \eqno(7.7)(a'')$$
We say that $w_\omega$ is a {\it three-complex
descent for ${\Cal L}'$}.  The eigenspaces of $T_{w_\omega}$ on the
span of ${\Cal L}$ and ${\Cal L}'$ are 
$$\text{$u^3$-eigenspace} = \langle {\Cal L} + {\Cal L}'\rangle, \quad
\text{$-1$-eigenspace} = \langle u^3{\Cal L} - {\Cal
L}'\rangle. \eqno(7.7)(b)$$ 

These two cases are the only possibilities when $\theta\alpha$ does
not belong to the span of $\alpha$ and $\beta$.  We consider next the
possibility $\theta\alpha$ {\it does} belong to this span. If
$\theta\alpha = (\alpha + \beta)$, applying $\theta$ 
gives $\alpha = \theta(\alpha+\beta)$, and therefore
$\theta\beta=-\beta$; so $\alpha$ is complex and $\beta$ is real,
contradicting the fact that $\alpha$ and $\beta$ must have the same
status. In the same way we rule out $\theta\alpha=
-(\alpha+\beta)$. Four possibilities remain: $\theta\alpha =
\pm\beta$, and $\theta\alpha = \pm\alpha$. We treat these next. 

We say that $w_\omega$ is a {\it three-semiimaginary ascent} if
$$\theta\alpha = \beta, \qquad \theta\beta = \alpha. \eqno (7.7)(c)$$
In this case of course $\theta(\alpha + \beta) = \alpha+\beta$, so
$\alpha+ \beta$ is an imaginary root; it turns out (by a simple
calculation in $SL(3)$) that it must be type II noncompact. (This case
arises in the example (9.2)(c) below, and the
calculations we omit can be done there.) The two
parameters $s\times {\Cal L}$ and $t\times {\Cal L}$ are interchanged by
$\sigma$. For $s\times{\Cal L}$, the simple root $\alpha+\beta$
corresponds to $t$, so the Cayley transform
$c_t(s\times {\Cal L})$ is double valued.  Evidently
$$\sigma\cdot c_t(s\times{\Cal L}) = c_s(t\times{\Cal L});$$
by calculation in $SL(3)$, we find that these two two-element sets
have a single element ${\Cal L}'$ in common. We define a {\it
single-valued Cayley transform}
$$c_{w_\omega}({\Cal L}) = c_t(s\times{\Cal L}) \cap c_s(t\times{\Cal
L})  = \{{\Cal L}'\}, \quad
{\Cal L}' \in {\frak D}^\sigma, \quad \ell({\Cal L}') = \ell({\Cal L}) + 2.
\eqno(7.7)(c')$$ 
A further calculation in $SL(3)$ leads to 
$$T_{w_\omega}({\Cal L}) = u{\Cal L} + (u+1){\Cal L}', \quad T_{w_\omega}({\Cal L}') = (u^3-u){\Cal L} + (u^3-u-1){\Cal L}'. \eqno(7.7)(c'')$$ 

We say that $w_\omega$ is a {\it three-real
descent for ${\Cal L}'$}; what this means is that
$\alpha$ and $\beta$ are real roots satisfying the parity condition
for ${\Cal L}'$.  We define a {\it
single-valued inverse Cayley transform}
$$c^{w_\omega}({\Cal L}') = s\times c^t({\Cal L}') = t\times c^s({\Cal
L}')  = \{{\Cal L}\} \in {\frak D}^\sigma. \eqno(7.7)(d)$$ 

The eigenspaces of $T_{w_\omega}$ on the 
span of ${\Cal L}$ and ${\Cal L}'$ are 
$$\text{$u^3$-eigenspace} = \langle {\Cal L} + {\Cal L}'\rangle, \quad
\text{$-1$-eigenspace} = \langle (u^2-u){\Cal L} - {\Cal
L}'\rangle. \eqno(7.7)(d')$$ 

We say that $w_\omega$ is a {\it three-imaginary noncompact ascent} if
$$\text{$\alpha$ and $\beta$ are noncompact imaginary
roots}. \eqno(7.7)(e)$$
Consequently $\alpha+\beta$ is compact imaginary. Because $W(K,H_i)$
must preserve the grading of the imaginary roots into compact and
noncompact, it follows that
$$W(K,H_i) \cap \langle s_\alpha,s_\beta\rangle =
\{1,s_{\alpha+\beta}\}.$$ 
In particular, the roots $\alpha$ and $\beta$ must be type I
noncompact imaginary. (This case arises in example $(9.2)(d)$ below.)
The Cayley transforms $c_s({\Cal L})$ and 
$c_t({\Cal L})$ are therefore single valued, and the two resulting
elements of ${\frak D}$ are interchanged by $\sigma$.  We define a
{\it single-valued Cayley transform}
$$c_{w_\omega}({\Cal L}) = s\times c_t({\Cal L}) = t\times c_s({\Cal
L}) = \{{\Cal L}'\} \in {\frak D}^\sigma, \quad \ell({\Cal L}') =
\ell({\Cal L}) + 2. \eqno(7.7)(e')$$ 
Calculation in $SL(3)$ (see $(9.2)(e)$ below) gives 
$$T_{w_\omega}({\Cal L}) = u {\Cal L} +(u+1){\Cal
L}',\quad T_{w_\omega}({\Cal L}') = (u^3-u){\Cal L} + (u^3-u-1){\Cal
L}' \eqno(7.7)(e'')$$ 
We say that $w_\omega$ is a {\it three-semireal descent for ${\Cal
L}'$}, and define a {\it single-valued inverse Cayley transform}
$$c^{w_\omega}({\Cal L}') = c^t(s\times{\Cal L}') \cap
c^s(t\times{\Cal L}') = \{{\Cal L}\}, \quad \ell({\Cal L}') =
\ell({\Cal L}) - 2. \eqno(7.7)(f')$$ 
The eigenspaces of $T_{w_\omega}$ on the span of ${\Cal L}$ and ${\Cal
L}'$ are  
$$\text{$u^3$-eigenspace} = \langle {\Cal L} + {\Cal L}'\rangle, \quad
\text{$-1$-eigenspace} = \langle (u^2-u){\Cal L} - {\Cal
L}'\rangle. \eqno(7.7)(f'')$$

We say that $w_\omega$ is a {\it three-real nonparity ascent for
${\Cal L}$} if $\alpha$ and $\beta$ are real non-parity roots. In this case
$$T_{w_\omega}({\Cal L}) =  -{\Cal L}, \quad \text{$u^3$-eigenspace} =
0, \quad \text{$-1$-eigenspace} = \langle{\Cal L}\rangle.\eqno(7.7)(g)$$

We say that $w_\omega$ is a {\it three-imaginary compact descent
for ${\Cal L}$} if $\alpha$ and $\beta$ are compact imaginary roots.
In this case
$$T_{w_\omega}({\Cal L}) = u^3{\Cal L}, \quad \text{$u^3$-eigenspace}
= \langle{\Cal L}\rangle, \quad \text{$-1$-eigenspace} = 0.\eqno(7.7)(h)$$

\head 8. Recursive algorithm\endhead
\subhead 8.1 \endsubhead
Recall from (5.2)(d) that the bar operator on $M$ satisfies
$$\DD({\cl}) = \sum_{\cl'\in\fD^\s;\cl'\preceq\cl}
\r_{\cl',\cl}{\cl'}, \qquad \r_{\cl,\cl} =
u^{-\dim[\cl]}. \eqno(8.1)(a)$$ 
By $(4.8)(e)$ we have also (writing $m=\ul(w_\omega)$)
$$u^{-m}(T_{w_\omega}+1)(\DD(\xi))=\DD((T_{w_\omega}+1)(\xi))
\qquad (\xi \in M). \eqno(8.1)(b)$$
In this section we will explain how to use these facts together with
the formulas of \S 7 to calculate the bar operator.

\proclaim{Proposition 8.2} In the setting (7.2), suppose $\cl \in
\fD^\s$. There are six mutually exclusive possibilities.

(a) {\it Good proper ascent.} The dimension of the support of
$(T_{w_\omega} + 1)(\cl)$ strictly exceeds that of $[\cl]$, and
involves one or two terms $\{\cl'_1,\cl'_2\}$ on a larger orbit, and
{\it no} other terms on an orbit of the same dimension. In any element of the
kernel of $T_{w_\omega}+1$, $\cl$ must appear as a multiple of
$$a\cl + b'_1\cl'_1 + b'_2{\Cal L}'_2,$$
with $a$, $b'_1$, and (if it is present) $b'_2$ all non zero. Therefore
(still in an element of the kernel of $T_{w_\omega}+1$)
$$\text{coefficient of $\cl$} = (\text{coefficient of $\cl'_1$})\cdot a/b'_1.$$

(b) {\it Bad proper ascent.} The dimension of the support of
$(T_{w_\omega} + 1)(\cl)$ strictly exceeds that of $[\cl]$, and
includes another term $\cl_2$ with $\dim[\cl_2] = \dim[\cl]$. 

(c) {\it Good proper descent.} The expression of
$(T_{w_\omega}+1)(\cl)$ involves only the term $\cl$ on the
orbit $[\cl]$, and one or two terms $\{\cl_1,\cl_2\}$ on
orbits strictly smaller than $[\cl]$.  In any element of the
image of $T_{w_\omega}+1$, $\cl$ must appear as a multiple of
$$\cl + \epsilon_1\cl_1 + \epsilon_2\cl_2,$$
with $\epsilon_1$ (and $\epsilon_2$ if it is present) equal to $\pm
1$. In particular.
$$(T_{w_\omega}+1)\cl_1 = x(\cl + \epsilon_1\cl_1 +
\epsilon_2\cl_2),$$
with $x\ne 0$. 

(d) {\it Bad proper descent.} The expression of
$(T_{w_\omega}+1)(\cl)$ involves {\it two} terms supported on the
orbit $[\cl]$. 

(e) We have $(T_{w_\omega} + 1)(\cl) = 0$. This is the case
of a real nonparity ascent, or $(7.5)(j''')$: type II imaginary with
$c_s({\Cal L})$ not fixed by $\sigma$.

(f) We have $(T_{w_\omega} + 1)(\cl) = (u^m + 1)\cl$. This is the
case of an imaginary compact descent, or $(7.5)(h''')$: type I real
with $c^s({\Cal L})$ not fixed by $\sigma$. 
\endproclaim
This is a summary of the calculations in \S 7. In type two, for
example, Proposition 8.2$(a)$ corresponds to the ascents in
$(7.6)(a)$, $(7.6)(c)$, $(7.6)(e)$, and $(7.6)(i)$; Proposition
8.2$(b)$ corresponds to the case $(7.6)(g)$; Proposition 8.2$(c)$
corresponds to the descents in $(7.6)(a)$, $(7.6)(c)$, (7.6)(g), and
(7.6)(i); and Proposition 8.2$(d)$ corresponds to $(7.6)(e)$.

\subhead 8.3 Reduction to a Levi subgroup \endsubhead
Before we formulate the recursion algorithm, it is useful to describe
a setting in which calculations (of the operator $\DD$ on $M$, and of
the polynomials $P^\s_{\cl',\cl}$) can be reduced to a Levi subgroup of
$G$.  Assume therefore that 
$$Q = LU \eqno(8.3)(a)$$
is a $(\sigma,\theta)$-stable Levi decomposition of a
$(\sigma,\theta)$-stable parabolic subgroup of $G$. We write
$\sigma_L$ for the restriction of $\sigma$ to $L$. Then there is a
natural $Q\cap K$-equivariant inclusion of flag varieties
$${\Cal B}_L \hookrightarrow {\Cal B}_G = {\Cal B}, \qquad B_L \mapsto
B_LU. \eqno(8.3)(b)$$
This inclusion gives rise to an isomorphism
$$K\times_{Q\cap K} {\Cal B}_L \simeq K\cdot {\Cal B}_L,
\eqno(8.3)(c)$$
and to inclusions
$${\frak D}_L \hookrightarrow {\frak D}, \qquad {\frak D}^{\sigma_L}_L
\hookrightarrow {\frak D}^\sigma.\eqno(8.3)(d)$$
Because $Q\cap K$ is parabolic in $K$, the variety $K\cdot {\Cal B}_L$
is a smooth $K$-stable closed subvariety of ${\Cal B}$.

\proclaim{Proposition 8.4} In the setting (8.3), complexes of $Q\cap
K$-equivariant constructible $\bbq$-sheaves on ${\Cal B}_L$ may be
identified with complexes of $K$-equivariant constructible
$\bbq$-sheaves on $K\cdot {\Cal B}_L$. This identification
respects Verdier duality, carries intersection cohomology sheaves to
intersection cohomology sheaves, and (in the language of disconnected
groups from (0.1)) carries $\widehat{L\cap K}$-equivariant sheaves to
$\hat K$-equivariant sheaves. In particular, if $\cl\in {\frak
D}_L^{\sigma_L}$, then (writing a superscript $G$ for the group)
$$P^{\sigma,G}_{\cl',\cl} = \cases P^{\sigma_L,L}_{\cl',\cl}, & \cl' \in
{\frak D}^{\sigma_L}_L \\
0, & \cl' \notin {\frak D}^{\sigma_L}_L . \endcases$$
$$\r^G_{\cl',\cl} = \cases \rho^{L}_{\cl',\cl}, & \cl' \in
{\frak D}^{\sigma_L}_L \\
0, & \cl' \notin {\frak D}^{\sigma_L}_L . \endcases$$
\endproclaim
All of this is an immediate consequence of the isomorphism (8.3)(c). 

On the level of representation theory in characteristic zero,
Proposition 8.4 corresponds to the cohomological induction functors
introduced by Zuckerman (see \cite{V}) to relate $({\frak l},L\cap
K)$-modules to $({\frak g},K)$-modules.

\subhead 8.5 How the recursion works \endsubhead
Calculating the bar operator on $M$ means calculating all of the
coefficients $\r_{\cl',\cl}\in \ca$.  We will explain how to do this by
(upward) induction on $\cl$; and then, for fixed $\cl$, by {\it downward}
induction on $\cl'$.

So fix $\cl$ and $\cl'$. We assume that we know how to calculate
$$\DD(\cl_1), \qquad \dim[\cl_1] < \dim [\cl]; \eqno(8.5)(a)$$
and that we know all the coefficients 
$$\r_{\cl'_1,\cl}, \qquad \dim[\cl'_1] > \dim[\cl'].\eqno(8.5)(a')$$
The first step in the recursion is to look for a proper good descent
$w_\omega$ for $\cl$ (Proposition 8.2(c)). If it exists, then we can
find one or two elements $\cl_j$ (with $[\cl_i]$ properly contained in the closure of $[\cl]$), signs $\epsilon_j$, and a nonzero
$x\in \ca$, so that
$$x \cl = (T_{w_\omega}+1)\cl_1 -
x(\epsilon_1\cl_1 + \epsilon_2\cl_2). \eqno(8.5)(b)$$
Applying $\DD$ to this formula yields 
$$\overline{x} \DD(\cl) = u^{-m}(T_{w_\omega}+1)\DD(\cl_1) -
\overline{x}\DD(\epsilon_1\cl_1 + \epsilon_2\cl_2). \eqno(8.5)(b')$$
This latter formula computes $\overline{x}\r_{\cl',\cl}$ in terms of various
$\r_{\cl'',\cl_i}$ (one needs all $\cl''$ so that $\cl'$ appears in
$(T_{w_\omega} + 1)\cl'$).

Recall that we write $B_j$ for a Borel subgroup in $[\cl]$. 
We may assume from now on that 
$$\text{there is no proper good descent for $\cl$.} \eqno(8.5)(c)$$
According to \cite{V}, Lemma 8.6.1, it follows that 
$$\gathered \text{the parabolic $Q=LU$ generated by}\\ \text{simple
real roots for $\cl$ is $\theta$-stable.}\endgathered\eqno(8.5)(c')$$
(In fact the absence of good descents means that all these simple real
roots must be either nonparity, or one-real type II, or one-real type
I with $c^s({\Cal L})$ not fixed by $\sigma$ (cf. $(7.5)(h''')$) or
two-real type II-II. But we do not yet need this more precise
information.) Proposition 8.4 reduces the calculation of $\DD$ to the
Levi subgroup $L$. We may therefore assume
$$\text{every simple generator $w_\omega$ is real for
$\cl$.}\eqno(8.5)(c'')$$

The next step in the recursion is to look for a proper good ascent
$w_\omega$ for $\cl'$ (Proposition 8.2(a)), with $\{\cl'_1,\cl'_2\}$
the corresponding terms on strictly larger orbits. If $w_\omega$ is a descent
for $\cl$, with $c^{w_\omega}(\cl) = \{\cl_1,\cl_2 \}$, then an easy
argument parallel to $(8.5)(b)$ computes $\r_{\cl',\cl}$ in terms of
$\DD(\cl_1)$, $\DD(\cl_2)$, and $\r_{\cl'_j,\cl}$; we omit the details.

Suppose therefore that
$$\text{$w_\omega$ is a proper good ascent for $\cl'$, and real
nonparity for $\cl$.} \eqno(8.5)(d).$$
According to \S 7, this means that $\cl$ belongs to the kernel of the
operator $T_{w_\omega} + 1$. By $(8.1)(b)$,
$$\DD\cl \in \ker (T_{w_\omega} + 1). \eqno(8.5)(d')$$
By Proposition 8.2(a), it follows that the appearance of $\cl'$ in
$\DD\cl$ must be as a multiple of $a\cl' + b'_1\cl'_1 + b'_2\cl'_2$.
Consequently 
$$\r_{\cl',\cl} = \r_{\cl'_1,\cl}\cdot a/b'_1, \eqno(8.5)(d'')$$
which is a recursion formula of the kind we want.  

It is immediate from $(8.5)(d')$ that 
$$\text{if $w_\omega$ is a compact descent for $\cl'$, then
$\r_{\cl',\cl} = 0$}.$$

Finally, we may assume (always in the setting $(8.5)(c)$) that
$$\gathered\text{there is no proper good ascent for $\cl'$, and}\\
\text{any compact descent for $\cl'$ is {\it not} real nonparity for
$\cl$.}\endgathered \eqno(8.5)(e)$$
According to \cite{V}, Lemma 8.6.1 (applied to $-\theta$), the first
condition implies that 
$$\gathered \text{the parabolic $Q'=L'U'$ generated by}\\ \text{simple imaginary roots for $\cl'$ is sent to its opposite by $\theta$.} \endgathered
\eqno(8.5)(e')$$ 
At this point we are able to work inside the ($\theta$-stable) Levi
subgroup $L'$; so we assume $L'=G$.  Our two parameters $\cl$ and
$\cl'$ satisfy 
$$\gathered \text{$[\cl]$ open, simple roots nonparity or
type II (see (7.5));}\\ \text{$[\cl']$ closed, 
simple roots compact or noncompact type I;} \\
\text{compact simple for $\cl'$} 
\longleftrightarrow \text{parity simple for $\cl$; and}\\
\text{noncompact simple for $\cl'$} 
\longleftrightarrow \text{nonparity simple for $\cl$.}\endgathered
\eqno(8.5)(e'')$$ 

Now the recursion is finished by
\proclaim{Proposition 8.6} In the setting (8.5), suppose that $\cl'
\ne \cl$. Then we can find
generators $w_{\omega_1}$ and $w_{\omega_2}$ so that 

(a) $w_{\omega_1}$ is a proper (imaginary) bad ascent for $\cl'$, and real
nonparity for $\cl$;

(b) $w_{\omega_2}$ is a compact descent for  $\cl'$, and a real
descent for $\cl$; and

(c) $w_{\omega_2}$ is an imaginary ascent for $w_{\omega_1} \times
\cl'$. 
\endproclaim

This is a purely structural assertion about involutive automorphisms
$\theta$ admitting both a maximal torus $H_s$ on which $\theta$ acts
by inversion, and one $H_c$ on which $\theta$ acts trivially.  We omit
the proof. 

Using the bad ascent $w_{\omega_1}$, the method of $(8.5)(d)$ gives a recursive
formula for $\r_{\cl',\cl} + \r_{w_{\omega_1}\times \cl',\cl}$.
The omitted argument before $(8.5)(d)$ gives (using $w_{\omega_2}$) a
recursive formula for $\r_{w_{\omega_1}\times \cl',\cl}$. Subtracting
these two formulas computes $\r_{\cl',\cl}$, and completes the algorithm.

\head 9. Examples\endhead
\subhead 9.1. Type $A_1\times A_1$\endsubhead Fix a fourth root
of one $i\in \kk$. Let 
$$\widetilde G = SL(2) \times SL(2),\qquad \widetilde\sigma(x,y) = (y,x)\quad
((x,y) \in \widetilde G), \eqno(9.1)(a)$$
$$\widetilde\theta\colon \widetilde G \rightarrow \widetilde G,\qquad
\widetilde\theta = \Ad\left( \pmatrix i & 0 \\ 0 & -i\endpmatrix,
\pmatrix i & 0 \\ 0 & -i\endpmatrix \right) \eqno (9.1)(b)$$
Then $\widetilde\sigma$ and $\widetilde\theta$ are commuting
involutive automorphisms of $\widetilde G$.  The group of fixed points
of $\widetilde\theta$ is
$$\widetilde K = \{\text{diagonal matrices in $\widetilde
G$\}}. \eqno(9.1)(c)$$ 
We will need two additional elements
$$ s_{\ell} = \left( \pmatrix 0 & 1 \\ -1 & 0\endpmatrix,
I\right),\quad  s_r = \left( I, \pmatrix 0 & 1 \\ -1 & 0\endpmatrix
\right). \eqno (9.1)(d)$$
The center of $\widetilde
G$ is naturally 
$$Z(\widetilde G) = \{\pm I\} \times \{\pm I\}; \eqno(9.1)(e)$$
$\widetilde\sigma$ acts by permuting the factors, and
$\widetilde\theta$ acts trivially.  There are three
$\widetilde\sigma$-stable subgroups of $Z(G)$:
$$Z_1 = \{(I,I)\}, \quad Z_2 = \{\pm(I,I)\}, \quad Z_3 = Z(G).$$
We can therefore make three examples of symmetric spaces with automorphisms
by dividing by subgroups of $Z(G)$: 
$$G_1 = \widetilde G, \qquad K_1 = \widetilde K \eqno(9.1)(f)$$
$$G_2 = \widetilde G/Z_2 , \qquad K_2 = \langle \widetilde
K, s_\ell s_r \rangle/Z_2 \eqno(9.1)(g)$$ 
$$G_3 = \widetilde G/Z_3 , \qquad K_3 = \langle \widetilde
K, s_\ell, s_r \rangle/Z_3 \eqno(9.1)(h)$$
The corresponding real semisimple groups are $Spin(3,1)$, $SO(3,1)$,
and $PSO(3,1)$.

Write $\phi_1$ for the standard Frobenius automorphism of $SL(2)$ over
$\FF_q$, raising each matrix entry to the $q$th power. The
quasisplit form of $G_1$ over $\FF_q$ is 
$$\align G_1(\FF_q) &= \{(x,y) \mid (x,y) =
(\phi_1(y),\phi_1(x)) \} \\
&= \{(x,\phi_1(x))\mid \phi_1^2(x) = x \} \simeq SL(2,\FF_{q^2}). \endalign$$ 
As usual we therefore have
$${\Cal B}(\FF_q) \simeq {\Bbb P}^1(\FF_{q^2}) \simeq \FF_{q^2} \cup \{\infty\}.$$
Clearly
$$K_1(\FF_q) \simeq \text{diagonal torus} \simeq \FF_{q^2}^\times.$$
The action of $z\in \FF_{q^2}^\times$ on the flag variety is by
multiplication by $z^2$. Therefore there are exactly four orbits of
$K_1(\FF_q)$ on ${\Cal B}(\FF_q)$:
$$\{0\}, \quad \{\infty\},\quad \{ \text{squares in
$\FF_{q^2}^\times$} \}, \quad  \{ \text{nonsquares in
$\FF_{q^2}^\times$} \},$$ 
That is, the space ${\Cal F}_s$ of $(6.1)(a)$ is four-dimensional. Using
the orbit description above, we can use the method of $(6.2)$ and $(6.3)$ to
calculate the action of ${\bold H}$, arriving at the formulas in
$(7.6)(e')$ and $(7.6)(f'')$.  A similar analysis applies to the quotient
groups $G_3$ and $G_2$, leading to the formulas of $(7.6)(g'')$,
$(h')$, $(i'')$, and $(j'')$. 

\subhead 9.2. Type $A_2$\endsubhead 
Define 
$$R = \pmatrix 0 & 0 & -1\\ 0 & 1 & 0 \\ -1 & 0 & 0 \endpmatrix, \quad
D = \pmatrix -1 & 0 & 0 \\ 0 & 1 & 0 \\ 0 & 0 & -1 \endpmatrix
\eqno(9.2)(a)$$
$$G = SL(3),\qquad \sigma(g) = R [{}^t g^{-1}] R^{-1}\quad
(g \in G). \eqno(9.2)(b)$$
The involutive automorphism $\sigma$ preserves the diagonal maximal
torus and the upper triangular Borel subgroup.  It may be written also
as 
$$\sigma\pmatrix a&b&c \\ d&e&f \\ g&h&i \endpmatrix
= \pmatrix ae-bd & af-cd & bf-ec \\ ah-bg & ai-cg & bi-ch \\
  dh - eg & di - gf & ei - fh \endpmatrix. \eqno(9.2)(b')$$
The right side is the action of $SL(3)$ on the three-dimensional
space $\bigwedge^2 \kk^3$ endowed with the basis $\{e_1\wedge e_2,
e_1\wedge e_3, e_2\wedge e_3\}$.
There are two interesting symmetric space automorphisms commuting with
$\sigma$. First is
$$\theta_s(g) = R[{}^t g^{-1}]R^{-1},\quad G^{\theta_s} = K_s = SO(3),
\eqno(9.2)(c)$$ 
the special orthogonal group with respect to the (maximally isotropic)
quadratic form 
$$\langle (a,b,c),(a',b',c') \rangle_O = -ac' -
ca' + bb'$$
The corresponding real semisimple group is $SL(3,{\RR})$. Second is
$$\theta_c(g) = D g D^{-1},\quad  G^D= K_c \simeq S(GL(2) \times
GL(1)). \eqno(9.2)(d)$$
(Here the $GL(2)$ block corresponds to the four corners of the matrix,
and the $GL(1)$ to the central element of the matrix. The remaining
four entries---$b$, $d$, $f$, and $h$, in the notation of
$(9.2)(b')$---are zero.)  The corresponding real semisimple group is
$SU(2,1)$.

The group $G(\FF_q)$ is $SU(3,\FF_q)$, the group of unitary
matrices of determinant one for the Hermitian form
$$\langle (a,b,c),(a',b',c') \rangle_U = -a \overline{c'} -
c\overline{a'} + b\overline{b'} \eqno(9.2)(e)$$
on $\FF_{q^2}^3$.  (Here we write $\overline a$ for the action of
the nontrivial Galois element in the quadratic extension of $\FF_q$ by $\FF_{q^2}$.) A maximally split torus is 
$$H(\FF_q) = \left\{\pmatrix a & 0 &0 \\ 0& a^{-1}\overline{a} &
0 \\ 0 & 0 & \overline{a}^{-1} \endpmatrix \mid a \in \FF_{q^2}^\times
\right\};$$ 
the upper triangular matrices in $SU(3,\FF_{q^2})$ comprise a
Borel subgroup. It follows easily that we may identify
$${\Cal B}(\FF_q) \simeq \{\text{isotropic lines in $\FF_{q^2}^3$}\}; \eqno(9.2)(e')$$
there are $q + 1 + (q^2-1)q = q^3 +1$ such lines, with
representative elements
$$\{(1,0,t)\mid \overline t = -t\}, \quad (0,0,1), \quad \{(x,1, y)
\mid x, y\in \FF_{q^2},\  x\overline y + y\overline x = 1
\}. \eqno(9.2)(e'')$$
(The counting is accomplished by noticing that there are $q$
possibilities for $t$, $q^2-1$ for $x$, and (for each fixed $x$) $q$
for $y$.)

The subgroup $K_c(\FF_q)$ is $S(U(2,\FF_q) \times U(1,\FF_q))$, with the first factor acting on the first and last
coordinates. It is therefore clear (by Witt's theorem that a unitary
group acts transitively on nonzero vectors of a given length) that
there are exactly two orbits 
of $K_c(\FF_q)$ on ${\Cal B}(\FF_q)$: the first $q+1$
isotropic lines (which live on the first and last coordinates), and
the remaining $q^3-q$ lines. These two orbits correspond by $(6.1)(a)$
to the local systems ${\Cal L}$ and ${\Cal L}'$ of $(7.7)(e')$, and
the formulas for the orbit cardinalities lead immediately (by $(6.2)$) to
$(7.7)(e'')$ and $(7.7)(f'')$. 

The subgroup $K_s(\FF_q)$ is $SO(3,\FF_q)$, the subgroup of
matrices with entries in $\FF_q$.  We omit the computation of the
orbits of $SO(3,\FF_q)$ on ${\Cal B}(\FF_q)$. It
turns out that there are exactly three: one of cardinality $q+1$, and
two of cardinality $(q^3-q)/2$.


\widestnumber\key{ALTV}
\Refs
\ref\key{\Un}\by J. Adams and M. van Leeuwen and P. Trapa and
D. Vogan  \paper Unitary representations of real reductive groups
\finalinfo{\tt arxiv:1212.2192}\endref
\ref\key{\BBD}\by A. Beilinson and J. Bernstein and P. Deligne \paper
Faisceaux pervers \jour Ast\`erisque \vol 100 \yr 1982 \endref
\ref\key{\IW}\by N. Iwahori\paper On the structure of the Hecke ring
of a Chevalley group over a finite field 
\jour J. Fac. Sci. Univ. Tokyo Sect. 1\vol10\yr 1964\pages 215--236\endref
\ref\key{\KL}\by D. Kazhdan and G. Lusztig\paper Schubert varieties
and Poincar\'e duality \paperinfo in {\it Geometry of the Laplace
Operator} \jour Proc. Symp. Pure  
Math.\vol36\publ Amer. Math. Soc.\yr1980\pages185--203 \eds
R. Osserman and A. Weinstein \endref 
\ref\key{\LC}\by G. Lusztig\paper Left cells in Weyl groups\inbook Lie groups representations, LNM 1024\publ
Springer Verlag\yr1983\pages99--111\endref
\ref\key{\CS}\by G. Lusztig\paper Character sheaves, V\jour
Adv. Math. \vol61 \yr1986 \pages103--155\endref
\ref\key{\QG}\by G. Lusztig\book Introduction to quantum
groups\bookinfo Progr. in Math. \vol110 \publ Birkh\"auser 
Boston \yr1993 \endref
\ref\key{\DG}\by G. Lusztig\paper On the representations of
disconnected reductive groups over $F_q$ \paperinfo in {\it Recent developments
in Lie algebras, groups, and representation theory}
\jour Proc. Symp. Pure Math.\vol86\publ
Amer. Math. Soc.\yr2012\pages227--246\eds K. Misra, D. Nakano, and
B. Parshall \endref 
\ref\key{\LV}\by G.Lusztig and D. Vogan\paper Singularities of closures of $K$-orbits on flag manifolds
\jour Invent. Math.\vol71\yr1983\pages 365--379\endref
\ref\key{\LVV}\by G. Lusztig and D. Vogan\paper Hecke algebras and
involutions in Weyl groups\jour
Bull. Inst. Math. Acad. Sin. (N.S.)\yr2012 \pages 323--354\vol 7\endref
\ref\key{\MA}\by H. Matsumoto\paper G\'en\'erateurs et relations des groupes de Weyl g\'en\'eralis\'es\jour C. R. Acad. Sci. Paris\vol258 \yr1964 \pages3419--3422\endref
\ref\key{V}\by D. Vogan\book Representations of reductive Lie
groups\bookinfo Progress in Mathematics\vol 15\publ Birkh\"auser\yr
1980 \endref 
\endRefs
\enddocument